\documentclass[12pt]{amsart}
\usepackage{amsmath}
\usepackage{amsxtra}
\usepackage{enumitem}
\usepackage{amscd}
\usepackage{amsthm}
\usepackage{hyperref}
\usepackage[T1]{fontenc}
\usepackage{amsfonts}
\usepackage{amssymb}
\usepackage{eucal}
\usepackage[all]{xypic}
\usepackage{color} 
\newtheorem{thm}[subsection]{Theorem}
\newtheorem{cor}[subsection]{Corollary}
\newtheorem{lem}[subsection]{Lemma}
\newtheorem{prop}[subsection]{Proposition}

\theoremstyle{definition}

\newtheorem{defn}[subsection]{Definition}

\newtheorem{rem}[subsection]{Remark}

\newtheorem*{thm*}{Theorem}
\newtheorem{question}[subsection]{Question}
\newtheorem{example}[subsection]{Example}

\newtheorem{problem}[subsection]{Problem}

\usepackage{yhmath}
\DeclareSymbolFont{largesymbols}{OMX}{yhex}{m}{n}
\DeclareMathAccent{\widetilde}{\mathord}{largesymbols}{"65}

\newcommand{\thmref}[1]{Theorem~\ref{#1}}
\newcommand{\secref}[1]{\S~\ref{#1}}
\newcommand{\lemref}[1]{Lemma~\ref{#1}}

\newcommand{\propref}[1]{Proposition~\ref{#1}}
\newcommand{\corref}[1]{Corollary~\ref{#1}}

\newcommand{\remref}[1]{Remark~\ref{#1}}
\newcommand{\exref}[1]{Example~\ref{#1}}

\newcommand{\nc}{\newcommand}
\nc{\renc}{\renewcommand}
\nc{\ssec}{\subsection}
\nc{\sssec}{\subsubsection} 
\nc\ol{\overline}
\nc\wt{\widetilde}
\nc\wh{\widehat}
\nc\wH{\widehat{H}}

\renc{\d}{{\delta}}
\nc{\Aa}{{\mathbb{A}}}
\nc{\Bb}{{\mathbb{B}}}
 \nc{\Gg}{{\mathbb{G}}}  
\nc{\Hh}{{\mathbb{H}}}
 \nc{\Nn}{{\mathbb{N}}}
\nc{\Rr}{{\mathbb{R}}}
\nc{\BV}{{\mathbb{V}}}
\nc{\BW}{{\mathbb{W}}}
\nc{\Zz}{{\mathbb{Z}}}
\nc{\Qq}{{\mathbb{Q}}}
\nc{\Ss}{{\mathbb{S}}}
\nc{\Cc}{{\mathbb{C}}}
\nc{\Ff}{{\mathbb{F}}}
 
\nc{\EL}{{L_\infty}}
 
\nc{\CA}{{\mathcal{A}}}
\nc{\CB}{{\mathcal{B}}}

\nc{\CE}{{\mathcal{E}}}
\nc{\CF}{{\mathcal{F}}}

 \def\tM{\widetilde{M}}

\nc{\Las}{\mathsf{Las}}
\nc{\CG}{{\mathcal{G}}}

\nc{\CL}{{\mathcal{L}}} 
\nc{\CC}{{\mathcal{C}}}
\nc{\CM}{{\mathcal{M}}}
\def\Mm{\CM}
\nc{\CN}{{\mathcal{N}}}
\nc{\Oog}{{\mathbb{O}}}
\nc{\Oo}{{\mathcal{O}}}
 \nc{\CQ}{{\mathcal{Q}}}
\nc{\CR}{{\mathcal{R}}} 
\nc{\CT}{{\mathcal{T}}}

\nc{\CV}{{\mathcal{V}}}
 
\nc{\CW}{{\mathcal{W}}}
\nc{\CZ}{{\mathcal{Z}}}

\nc{\cM}{{\check{\mathcal M}}{}}
\nc{\csM}{{\check{\mathcal A}}{}}
\nc{\oM}{{\overset{\circ}{\mathcal M}}{}}
\nc{\obM}{{\overset{\circ}{\mathbf M}}{}}
\nc{\oCA}{{\overset{\circ}{\mathcal A}}{}}
\nc{\obA}{{\overset{\circ}{\mathbf A}}{}}
\nc{\ooM}{{\overset{\circ}{M}}{}}
\nc{\osM}{{\overset{\circ}{\mathsf M}}{}}
\nc{\vM}{{\overset{\bullet}{\mathcal M}}{}}
\nc{\nM}{{\underset{\bullet}{\mathcal M}}{}}
\nc{\oD}{{\overset{\circ}{\mathcal D}}{}}
\nc{\obD}{{\overset{\circ}{\mathbf D}}{}}
\nc{\oA}{{\overset{\circ}{\mathbb A}}{}}
\nc{\op}{{\overset{\bullet}{\mathbf p}}{}}
\nc{\cp}{{\overset{\circ}{\mathbf p}}{}}
\nc{\oU}{{\overset{\bullet}{\mathcal U}}{}}
\nc{\oZ}{{\overset{\circ}{\mathcal Z}}{}}
\nc{\ofZ}{{\overset{\circ}{\mathfrak Z}}{}}
\nc{\oF}{{\overset{\circ}{\fF}}}

\nc{\fa}{{\mathfrak{a}}}
\nc{\fb}{{\mathfrak{b}}}
\nc{\fg}{{\mathfrak{g}}}
\nc{\fgt}{{\fg}_!}
\nc{\fgl}{{\mathfrak{gl}}}
\nc{\fh}{{\mathfrak{h}}}
\nc{\fj}{{\mathfrak{j}}}
\nc{\fm}{{\mathfrak{m}}}
\nc{\ft}{{\mathfrak{t}}}
\nc{\fn}{{\mathfrak{n}}}
\nc{\fu}{{\mathfrak{u}}}
\nc{\fp}{{\mathfrak{p}}}
\nc{\fr}{{\mathfrak{r}}}
\nc{\fs}{{\mathfrak{s}}}
\nc{\fsl}{{\mathfrak{sl}}}
\nc{\hsl}{{\widehat{\mathfrak{sl}}}}
\nc{\hgl}{{\widehat{\mathfrak{gl}}}}
\nc{\hg}{{\widehat{\mathfrak{g}}}}
\nc{\chg}{{\widehat{\mathfrak{g}}}{}^\vee}
\nc{\hn}{{\widehat{\mathfrak{n}}}}
\nc{\chn}{{\widehat{\mathfrak{n}}}{}^\vee}

\nc{\fA}{{\mathfrak{A}}}
\nc{\fB}{{\mathfrak{B}}}
\nc{\fD}{{\mathfrak{D}}}
\nc{\fE}{{\mathfrak{E}}}
\nc{\fF}{{\mathfrak{F}}}
\nc{\fG}{{\mathfrak{G}}}
\nc{\fK}{{\mathfrak{K}}}
\nc{\fL}{{\mathfrak{L}}}
\nc{\fM}{{\mathfrak{M}}}
\nc{\fN}{{\mathfrak{N}}}
\nc{\fP}{{\mathfrak{P}}}
\nc{\fU}{{\mathfrak{U}}}
\nc{\fV}{{\mathfrak{V}}}
\nc{\fZ}{{\mathfrak{Z}}}

\nc{\bb}{{\mathbf{b}}}
\nc{\bc}{{\mathbf{c}}}
\nc{\bd}{\partial}
\nc{\be}{{\mathbf{e}}}
\nc{\bj}{{\mathbf{j}}}
\nc{\bn}{{\mathbf{n}}}   \nc{\bm}{{\mathbf{m}}}
 \nc{\de}{{\bar{d}}}
\nc{\bp}{{\mathbf{p}}}
\nc{\bq}{{\mathbf{q}}}
\nc{\bF}{{\mathbf{F}}}
\nc{\bu}{{\mathbf{u}}}
\nc{\bv}{{\mathbf{v}}}
\nc{\bx}{{\mathbf{x}}}
\nc{\bs}{{\mathbf{s}}}
\nc{\by}{{\bar{y}}}
\nc{\bw}{{\mathbf{w}}}
\nc{\bA}{{\mathbf{A}}}
\nc{\bK}{{\mathbf{K}}}
\nc{\bI}{{\mathbf{I}}}
\nc{\bB}{{\mathbf{B}}}
\nc{\bG}{{\mathbf{G}}}
\nc{\bC}{{\mathbf{C}}}
\nc{\bD}{{\mathbf{D}}}
 
\nc{\bH}{{\mathbf{H}}}
\nc{\bM}{{\mathbf{M}}}
\nc{{\tN}}{{\widetilde{N}}}
\nc{\bV}{{\mathbf{V}}}
\nc{\bU}{{\mathbf{U}}}
\nc{\bL}{{\mathbf{L}}}
\nc{\bT}{{\mathbf{T}}}
\nc{\bW}{{\mathbf{W}}}
\nc{\bX}{{\mathbf{X}}}
\nc{\bY}{{\mathbf{Y}}}
\nc{\bZ}{{\mathbf{Z}}}
\nc{\bS}{{\mathbf{S}}}
\nc{\bSi}{{\bar{\Sigma}}}
\nc{\sA}{{\mathsf{A}}}
\nc{\sB}{{\mathsf{B}}}
\nc{\sC}{{\mathsf{C}}}
\nc{\sD}{{\mathsf{D}}}
\nc{\sF}{{\mathsf{F}}}
\nc{\sG}{{\mathsf{G}}}
\nc{\sK}{{\mathsf{K}}}
\nc{\sM}{{\mathsf{M}}}
\nc{\sO}{{\mathsf{O}}}
\nc{\sQ}{{\mathsf{Q}}}
\nc{\sP}{{\mathsf{P}}}
\nc{\sZ}{{\mathsf{Z}}}
\nc{\sfp}{{\mathsf{p}}}
\nc{\sr}{{\mathsf{r}}}
\nc{\sg}{{\mathsf{g}}}
\nc{\sff}{{\mathsf{f}}}
\nc{\sfb}{{\mathsf{b}}}
\nc{\sfc}{{\mathsf{c}}}
\nc{\sd}{{\ltimes}}

 \def\wB{\widehat{B}}

\nc{\Pp}{{\mathop{\operatorname{\rm P}}}}
  \nc{\vol}{{\mathop{\operatorname{\rm vol\,}}}}
\nc{\co}{{\mathop{\operatorname{\rm Core\,}}}}
  \nc{\gal}{{\mathop{\operatorname{\rm Gal\,}}}}
  \nc{\cl}{{\mathop{\operatorname{\rm cl}}}}
  \nc{\disc}{{\mathop{\operatorname{\rm disc}}}}
  \nc{\Sym}{{\mathop{\operatorname{\rm Sym}}}}
   \nc{\Aut}{{\mathop{\operatorname{\rm Aut}}}}
 \nc{\Spec}{{\mathop{\operatorname{\rm Spec}}}}
  \nc{\spec}{{\mathop{\operatorname{\rm Spec}}}}
\nc{\Ker}{{\mathop{\operatorname{\rm Ker}}}}
 \nc{\dom}{{\mathop{\operatorname{\rm dom}}}}
\nc{\End}{{\mathop{\operatorname{\rm End}}}}
 \nc{\Hom}{\operatorname{Hom}}
 \nc{\GL}{{\mathop{\operatorname{\rm GL}}}}
 \nc{\Id}{{\mathop{\operatorname{\rm Id}}}}
 \nc{\rk}{{\mathop{\operatorname{\rm rk}}}}
 \nc{\length}{{\mathop{\operatorname{\rm length}}}}
\nc{\supp}{{\mathop{\operatorname{\rm supp} \, }}}
\nc{\val}{{\rm val}}
 
\nc{\res}{{\mathop{\operatorname{\rm res}}}}

\def\Ind#1#2#3{{#1} {\downarrow}_{#3} {#2} }

\def\meet{\cap}
\def\union{\cup}

\def\si{\sigma}
\def\g{\gamma}
\def\G{\Gamma}

\def\m{\smallsetminus}

\nc{\seq}[1]{\stackrel{#1}{\sim}}

\def\inv{^{-1}}
\def\claim#1{\smallskip {\noindent \bf Claim #1.\quad}}

\def\beq#1{\begin{equation} \label{ #1}}
\def\eeq{\end{equation}}

\def\prf{\begin{proof}}
\def\pv{\end{proof} }
 \def\eprf{\end{proof} }

\def\acl{\mathop{\rm acl}\nolimits}
 \def\dcl{\mathop{\rm dcl}\nolimits}

 \renc{\b}{{\beta}}

\def\Ind#1#2{#1\setbox0=\hbox{$#1x$}\kern\wd0\hbox to 0pt{\hss$#1\mid$\hss}
\lower.9\ht0\hbox to 0pt{\hss$#1\smile$\hss}\kern\wd0}

 \def\Om{\Omega}
 \def\om{\omega}

\def\two{\mathsf{\bf{2}}}

 \author{Jamshid Derakhshan}
 \address{St Hilda's College, University of Oxford, Cowley Place, Oxford OX4 1DY, and St Peter's College, University of Oxford, New Inn Hall Street, Oxford OX1 2DL}
 \email{derakhsh@maths.ox.ac.uk}
 
 \author{Ehud Hrushovski}
 \address{Mathematical Institute, University of Oxford, Andrew Wiles Building, Radcliffe Observatory Quarter Woodstock Road, Oxford OX2 6GG}
 
\email{Ehud.Hrushovski@maths.ox.ac.uk}
   \dedicatory{To Katrin, with admiration and  friendship}
 
 \title{Imaginaries, products and the adele ring}  
\date{30 Jan., 2026}

 \def\afin{\Aa^{fin}}

  \def\tT{\widetilde{T}}

\setlist[itemize]{leftmargin=*}

\begin{document}

  \begin{abstract} We describe the imaginary sorts of infinite products in terms of imaginary sorts of the factors (\thmref{1}).  We extend the result to certain reduced powers (\thmref{2}) and  then to infinite products $\prod_{i\in I} M_i$ enriched with a predicate for the ideal of finite subsets of $I$ (\thmref{3}).   As a special case, using the Hils-Rideau-Kikuchi uniform $p$-adic elimination of imaginaries, we
 find the imaginary sorts of the ring of rational adeles (\exref{adeles}).     Our methods include the use of the Harrington-Kechris-Louveau Glimm-Efros dichotomy both for transitioning from monadic second order imaginaries to first-order reducts, and for proving a certain ``one-way'' model-theoretic orthogonality within the adelic imaginaries.
 
  \end{abstract}
 \maketitle

 \section{Introduction}

  \ssec{Products and Boolean-valued models}   The theory of a power $M^I$ of an algebraic structure $M$ was determined by Mostowski \cite{mostowski},
  both for full  and weak powers (such as direct powers and direct sums of a group.)    It was systematized and generalized to reduced products  by Feferman and Vaught \cite{FV}.   To understand such structures it is best to add to the picture  the Boolean algebra $B$ of all subsets of $I$
    (or all finite or cofinite subsets of $I$, in the case of weak products), 
and view $M^I$ as a $B$-valued model of the theory of $M$; see \secref{boolean}.   If the Boolean algebra is enriched with an ideal $J$, it becomes possible
to interpret restricted powers (see \secref{adeles}), weak products or reduced products.    In fact under weak conditions on $M$,  the Boolean algebra and attendant
ideal are  already interpretable in $M^I$, so it is more a matter of displaying than of really adding the Boolean apparatus.
   
An interesting example    is  
 the structure $(\Nn, \cdot)$ of multiplication on the integers; the fundamental theorem of arithmetic identifies it with  the weak power of Pressburger arithmetic 
 $(\Nn,+)$.    The study of  $(\Nn, \cdot)$, and a procedure for interpreting  
  Boolean algebra operations within it, go  back to \cite{leibniz}.    In modern terms,  $(\Nn, \cdot)$ interprets in the first-order sense the Boolean algebra of finite and cofinite subsets of $\Nn$, via square-free integers and divisibility.    Leibniz makes use not of this directly,    but of a closely related  ind-definable {\em free} Boolean algebra, generated by the primes.    The mathematical reason for this choice is clear:  his goal was to reduce all valid Aristotelian syllogisms
   to numerical calculations of multiplication and greatest common divisor operations on natural numbers; it is the free Boolean algebra on given unary predicates that captures  Aristotelian logic.     
 The first-order theory of $(\Nn,\cdot)$ was fully analyzed  by Skolem in 1930, following on his 1919 work on the  theory of 
 $\{0,1\}^\Nn$ - equivalently the theory of infinite atomic Boolean algebra.    (References to Skolem's work can be found in \cite{mostowski}.)

   The Boolean-valued approach fits into  a series of richer examples, including the Scott-Solovay-Shoenfield Boolean-valued approach to forcing, and the theory of the algebraic integers viewed as a Boolean-valued algebraically closed valued field.    It was  also a component of Keisler's theory of randomization \cite{keisler}, later transposed into
continuous logic by Ita\"i 
Ben Yaacov and Keisler.    In the present paper we treat only plain  Booleanizations, enriched at most by an ideal.   
 Guided by the case of the ring of adeles, we consider  the zero ideal or ideals $J$ such that the algebra is atomless modulo $J$.
  
Many thanks to Mariana Vicaria for her extended and comments on an earlier version of this paper,  to George Peterzil and to Kobi Peterzil, 
Mira Tartarotti,   and to the 
anonymous referee; the paper was  much improved by his or her suggestions.
  
  \ssec{Imaginaries} \label{imaginaries}
  An imaginary sort of a first-order structure $M$ is a possible universe of a structure interpretable in $M$.  It can be presented as a quotient $D/E$,
  with $D$ a definable set (perhaps of tuples) and $E$ a definable equivalence relation on $D$.  
    The significance of imaginary sorts was   established in the abstract setting of Classification Theory by Shelah \cite{shelah}; subsequently all notions of stability theory, especially algebraic closure, were interpreted 
   ``in $T^{eq}$'', i.e. with imaginary sorts taken into account.     They were 
     first studied for specific theories by Poizat \cite{poizat}; he showed notably that the theories of algebraically closed fields and of 
     differentially closed fields {\em eliminate imaginaries}, i.e. each imaginary sort is in definable bijection with    a definable set of $n$-tuples.  
     Poizat also  
  pointed  out that elimination of imaginaries for algebraically closed fields encapsulates, modulo general principles, classical Galois theory for fields.  
  
  Other theories have imaginary sorts that differ qualitatively from definable sets; an example is the value group of a valued field such as $\Qq_p$.  
 In this case it is a question of finding a useful {\em basis} for the imaginary sorts, meaning every imaginary sort can be definably embedded into a product of the basis elements.   In the case of $\Qq_p$, these would be the {\em geometric sorts} $GL_n(\Qq_p)/GL_n(\Zz_p)$
  of lattices in $\Qq_p^n$.
  
  The theory of an infinite set $V$, with no relations other than equality, does not admit elimination of imaginaries; since one has the imaginary sorts
  $V^n /H$, for any subgroup $H$ of the symmetric group $Sym(n)$; they code finite subsets of $V^n$.   
    However these are the only imaginary sorts one needs.    The same situation occurs in theories of  vector spaces over a given field.
 
  In general, a theory with universe $V$ has {\em weak elimination of imaginaries} if the sorts $V^n/H$  form a basis for the imaginaries.

 At this point we can give a formulation of  our main theorems.   Let $T$ be a first-order theory in a language $L$; $T$ may be incomplete.  We assume for this introduction that
 each sort $S$  has at least one $0$-definable element in $T$ (assumption $\clubsuit$; 
see \remref{clubrem}). 
  Let $I$ be an infinite index set, 
 and let $(M_i: i \in I)$ be a family of models of $T$.   Let $M^{\sharp} = \prod_{i \in I} M_i$ with the natural structure
 (explained precisely in \exref{products}), including a sort $B$ for the  Boolean algebra $\Pp(I)$.
     For an equivalence relation $E$ on a finite product $D$ of sorts of $T$,   let $E^{\sharp}$ and $E^{\flat}$ be  equivalence relations on $D(M)$ defined by:  $xE^{\sharp}y$ iff
     $x(i) E y(i)$ for all $i \in I$, and $x E^\flat y$ iff $x(i) E y(i) $ for   all but finitely many $i \in I$.   Let $M^{\flat}$ be the expansion of $M$ by 
  a predicate for   the ideal $I_{fin} \subset B$ of finite subsets of $I$.     
 \begin{thm*}  \item   $M^{\sharp}$ admits weak EI  to the sorts $D/E^{\sharp}$ (for any $D,E$ as above.)
 \item   $M^{\flat}$ admits weak EI  to the sorts $D/E^{\sharp}$ and $D/E^{\flat}$ (for any $D,E$ as above.) 
 \end{thm*}
 
 This will be proved as Theorems \ref{1} and   \ref{3} (applied to $T^{eq})$.     

\begin{rem} \label{clubrem}
 In case assumption $\clubsuit$ fails, \thmref{1} still implies
 weak elimination of imaginaries once 
 one adds a sort for the partial product $\Pi_{i \in u} M_i$, uniformly in elements $u$ of the Boolean algebra $B$. More precisely, it suffices
 to add a sort for  $\Pi_{i \in u} S(M_i)$  
for each sort $S$ of $L$ without a 0-definable element in $T$.  
 
 This follows from the present
 statement of \thmref{1}, applied to a theory bi-interpretable with $T$ where one adds a new "ideal element" $*_S$ of each such sort $S$.   An element $(e_i)$ of the full product $\Pi_i S(M_i) \union *_S$ can be identified
 with an element of $\Pi_{i \in u} S(M_i)$ obtained by deleting all ideal elements among the $e_i$ along with their indices.

 Similar remarks are valid for Theorems \ref{2} and \ref{3}.
  \end{rem}

    Weak elimination of imaginaries for Skolem arithmetic, $Th(\Nn, \cdot,1)$, was proved by Atticus Stonestrom 
  \cite{atticus}.  This can be viewed as a precursor and essentially a special case of  \thmref{3}, since  $(\Nn,\cdot,1)$   is interpretable in $(\Zz, P(I), I_{fin})$ where    $I$ is the set of primes, as the set of positive elements of $\Zz^P$ with  finite support.  
  
    An earlier result of Newelski and Wencel (\cite{NW},  Cor. 2.5)
  proves weak EI  for Boolean algebras with finitely many atoms.     Still earlier, as Wencel points out,  Truss \cite{truss} (Cor. 3.8) proved the small index property for  countable atomless Boolean algebras, in a form that includes  weak EI for this theory.     Wencel
     \cite{wencel} then proved the atomic case.

  In \secref{sec:atomic} we will prove our first general result reducing imaginaries in product structures to imaginaries of the individual structures.   The  fundamental case in our treatment
  is $\Om^\Nn$, where $\Om$ is a finite set with no additional structure.  While it may well be possible to handle this case explicitly, we handle it with tools
  of liaison groups originating in stability theory,  whose appearance and effectiveness in this second-order environment initially surprised us.  In particular, we give some structure theorems for interpretable groups over the  infinite atomic Boolean algebra.
  This can be compared to the determination of connected $\infty$-definable groups in the stable continuous logic theory of probability logic \cite{berenstein}.

  Along with  liaison groups, we will use tools from descriptive set theory.   A Borel equivalence relation $E$ on a Polish space $Y$ 
 is called {\em smooth} if there exists a 
  Borel map $f$ from $Y$ to a Polish space,  such that $f(x)=f(y)$ iff $xEy$.    In work of Becker and Kechris \cite{bk}, this is shown to give a very robust criterion for a Borel subgroup to be closed.    Roughly speaking, this will allow us to move from questions on imaginaries of a product $\prod_{i\in I} M_i$ to questions on 
  reducts of a co-product $\coprod_{i\in I} M_i$.   The latter remains stable in the critical case where the $M_i$ are finite, allowing the use of tools from geometric stability theory.

These descriptive-set theoretic ideas  will  also give us an important model-theoretic orthogonality principle.  
 Let $I_{fin}$ be the ideal of finite subsets of $I$, and $E_{fin}$ the associated equivalence relation:  
 $aE_{fin}  b$ iff $\{i\in I: a_i \neq b_i\} \in I_{fin}$. 
     In case $M_i$ is countable (but with at least two elements),  $E_{fin}  $ is a version of the paradigmatic non-smooth Borel equivalence   relation $E_0$;     Harrington, Kechris and Louveau  show in \cite{hkl} that a Borel equivalence relation is non-smooth if and only if 
  $E_{0}  $ can be reduced to it.  In the adelic case, where $M_i$ is a local field, or more generally an uncountable Borel structure, $E_{fin}$ gives the equivalence relation  known in descriptive set theory as $E_1$
  (which is even less smooth.)     On the other hand, we will see that the equivalence relations definable on a single $E_{fin}  $-class are all smooth, in these cases.
  It will follow from this that there are no nontrivial  definable maps  or finitely multivalued maps from the reduced power $M/E_{fin}  $
 into $M$ (\propref{function}); this conclusion will hold in any model.  This will allow us to study the imaginaries separately in each $E_{fin}  $- class  on various definable sets, equivalence relation, and in the quotient, and glue together the information.  

The elimination of imaginaries can presumably also be achieved in a more hands-on manner,   along the lines of \cite{wencel}, \cite{atticus}.
But our purpose here is to develop tools that may serve in more complex environments.   In \cite{bh}, the same descriptive set-theoretic methods
are used to show in far greater generality that smoothness coincides with an appropriately generalized notion of 2nd-order elimination of imaginaries.

  Our test case at present was the ring of rational adeles $\Bbb A_{\Qq}$, and the above general results lead   in particular to describe the imaginary sorts of the adele ring in \secref{adeles}.
  They consist of the uniform imaginary sorts for the $p$-adics, namely the ring itself and  the codes for finite sets of adelic lattices in $\Aa_{\Qq}^n$, but each appearing also in a second blurred version,    identifying two sequences that agree for almost all primes.    Thanks to work of Hils and Rideau-Kikuchi, a uniform description of the imaginaries of $\Qq_p$ in a natural language,
  roughly speaking the field language along with certain canonical generators for certain tame Galois groups.   We thus obtain
  a description of the adelic imaginaries in the same language.  
  
    It would be very interesting to find a basis for the imaginaries for  the algebraic integers, where a `discrete' field is also endowed from structure
    arising from embedding into an adele ring (see Problem \ref{algint}).

 \section{Preliminaries}
  
 Let ${L}$ be a language, possibly many-sorted.

 Let ${T}$ be a   first-order theory of ${L}$.  It will be convenient to assume (and can easily be arranged by adding sorts)
 that  any finite product of sorts stands in 0-definable bijection with a sort. Let $L_{BA}=\{\cup,\cap,0,1,-\}$ be the language of Boolean algebras. We denote the Boolean complement of an element $a$ by $1-a$. 
 
 Recall that in a given structure $M$, the {\em algebraic closure} of a subset $A$ is
 the union of all finite $M_A$-definable sets; while the {\em definable closure} is the set of all $M_A$-definable elements, equivalently
 the union of all $M_A$-definable one-element sets.   
 We denote algebraic and definable closure by  $\acl$ and $\dcl$, respectively.   While our goal is to eliminate abstract imaginary sorts
 in favour of concrete ones, as noted in \secref{imaginaries}, it will sometimes be useful to consider (abstract) imaginary elements; the notions of $\acl,\dcl$ apply also to $M^{eq}$.   Note that  
the relation $a \in \acl(b)$ between given elements is unambiguous - in case $a,b$ are `real', the relation holds in $M$ if and only if it holds
in $M^{eq}$,    and likewise for $\dcl$.    In one or two arguments it is essential to take the closure in $M$ - restrict to `real' elements;   we will then explicitly  say so.
 
\ssec{Language of $T^{boole}$}  \label{boolean}

Let $L^{boole}$ be a language with the sorts  of $L$ and one additional sort $B$ equipped with the language $L_{BA}$.  
(For simplicity we will use notation as if  $L$ has a single sort $K$.)   On $B$ we have the operations of a Boolean algebra,
and the equality relation.    In addition,
  for any   formula $\phi(x_1,\ldots,x_n)$ of $L$, we have a function symbol 
$[\phi]: K^n \to B$ with the same variables.  We view $[\phi]$ as giving the Boolean truth value of $\phi$. 
In (2-4) below, all variables $x,y,a,b,\ldots$  range over the sorts of $L$.

\ssec{Axioms of $T^{boole}$}  \label{axioms1}   
\begin{enumerate}
\item   $B$ is a Boolean algebra.  
\item  If $T \models \forall x \phi$, then for all $a$ from $K$, $[\phi](a)=1$; where $x=(x_1,..,x_n)$ and $a=(a_1,...,a_n)$.
\item For any tuple $a$ from $K$, $\phi \mapsto [\phi](a)$ is a Boolean homomorphism.   Also $[,]$ respects dummy variables, i.e. if $\phi(x,y):=\psi(x)$ 
then $[\phi](a,b)=[\psi](b)$. (We will sometimes denote $[\phi](z)$ by $[\phi(z)]$ when this simplifies the notation.)   
\item   \label{lg1}  Assume $T \vdash (\forall y)( \psi_i(y) \implies (\exists x) \phi_i(x,y))$, $i=1,\ldots,n$.  Then 
\[(\forall b_1,\ldots,b_n)( \bigwedge_{i \neq j} b_i \meet b_j = 0 \wedge \bigwedge_{i\leq n} [ \psi_i(y) ] \geq b_i \implies
 (\exists x) \bigwedge_{i\leq n} [\phi_i(x,y)] \geq b_i  ) \] 
 
\end{enumerate}

 Axiom (4) can also be presented in two parts:
 
 (4a) [Glueing]   Whenever $b_1,\ldots,b_n$ form a partition of unity in $B$, and $a_1,\ldots,a_n\in K$ are elements, there exists an element $a$
 with $[a=a_i] \geq b_i$.  
 
 (4b)   If $T \models (\forall x)(\exists y) \phi(x,y)$ then $T^{boole}$ includes the axiom $(\forall x)(\exists y) [\phi(x,y)] =1$.

 It follows from these axioms that in fact $T^{boole} \models [\psi]=1$ whenever $T \models \psi$ (for $\psi$ an $L$-sentence).
 
 We call a sequence of $L$-formulas $\phi_1(x),\dots,\phi_k(x)$ {\it a partition of unity in $T$} if for $i \neq j$, 
 $T \models \neg (\exists x)(\phi_i \wedge \phi_j)$, and  $T \models(\forall x)( \bigvee_{i\leq k} \phi_i(x) )$

  The theorem below can be viewed as a slightly generalized form of the Feferman-Vaught theorem \cite{FV}.  In \cite{keisler},   this is taken further by 
  placing a measure on the Boolean algebra (we do not do so here.)

   \begin{thm}\label{FV}   Let $\alpha(x,u)$  be a formula of $T^{boole}$,
  where $x$ is a tuple of variables of   sorts of $T$, and   $u$ is a tuple of variables of sort $B$.   Then 
  for some $L$-formulas $\phi_1(x),\dots,\phi_k(x)$, that form a partition of unity in $T$, and some $L_{BA}$-formula $\beta$, 
  \[ T^{boole} \models \alpha \iff \beta( [\phi_1(x)],\ldots,[\phi_k(x)],u). \]
  $T^{boole}$ is  complete modulo  the theory of $B$ in the language of
Boolean algebras $L_{BA}$ with additional constants, namely the elements $[\psi]$ for sentences $\psi$ of $L$.

Consequently $B$  is stably embedded.   The induced structure is just the Boolean algebra structure along with the above constants.

All these statements remain true if $B$ is an expansion of a Boolean algebra.
\end{thm}
  \prf   
We start with 

\claim{}Given $m'\geq 1$, $L$-formulas $\phi'_1(x),\dots,\phi'_{m'}(x)$ and a Boolean formula $\beta'$ in $m'$ variables, there are $L$-formulas $\phi_1(x),\dots,\phi_{m}(x)$ (in the same free variables $x$) for some $m\geq 1$, and a Boolean formula $\beta$ in $m$ variables such that 
 \[ \beta([\phi_1(x)],\dots,[\phi_m(x)]) \iff \beta'([\phi'_1(x)]\,\dots,[\phi'_{m'}(x)])\] 
 and such that $(\phi_1(x),\dots,\phi_m(x))$ is a partition of unity in $T$. 
  
\proof Let $m=2^{m'}$ and $s_1,\dots,s_m$ be a list of subsets of $\{1,\dots,m'\}$. For $k\leq m$, let 
 \[\phi_k:=\bigwedge_{j\in s_k} \phi'_j \wedge \bigwedge_{j\in \{1,\dots,m'\}\setminus s_k} \neg \phi'_j.\]
 
 For $l\leq m'$, let $r_l:=\{k: k\leq m,~ l\in s_k\}$, and define the Boolean formula $\beta:=\beta'(\bigvee_{k\in r_1} u_k,\dots,\bigvee_{k\in r_{m'}} u_k)$. Then $(\phi_1,\dots,\phi_m)$ is a partition of unity in $T$ and $\beta([\phi_1(x)],\dots,[\phi_m(x)])$ and $ \beta'([\phi'_1(x)]\,\dots,[\phi'_{m'}(x)])$ are equivalent. This proves the Claim.
 
 We now give the argument for QE.  It suffices to eliminate quantifiers in existential formulas, i.e. given $\exists x \beta'([\phi'_1(x,z)],\dots,[\phi'_{m'}(x,z)])$, where $\phi'_1,\dots,\phi'_{m'}$ are $L$-formulas, $x$ and $z$ are tuples of $L$-variables, and $\beta'$ is a Boolean formula;  one can construct $L$-formulas 
 $\phi_1(z),\dots,\phi_m(z)$ for some $m\geq 1$ and an Boolean formula $\beta$ such that
 \[(*)  \ \ \ \ \ T^{boole}\models \exists x \beta'([\phi'_1(x,z)],\dots,[\phi'_{m'}(x,z)]) \iff \beta([\phi_1(z)],\dots,[\phi_m(z)]).\]
 To show this we may assume that $(\phi'_1,\dots,\phi'_{m'})$ is a partition of unity in $T$ by the above. Let $m=m'$ and 
 $\phi_j :=\exists x \phi'_j$ for $j\leq m$. Let $Part(u_1,\dots,u_m)$ express that $u_1,\dots,u_m$ form a partition of unity in the Boolean algebra. Let 
 \[\beta(u_1,\dots,u_m):=\exists w_1 \dots \exists w_m (Part(w_1,\dots,w_m) \wedge \bigwedge_{j\leq m} (w_j \leq u_j) \wedge \beta'(w_1,\dots,w_m)).\]
 
 Suppose $\mathcal{M}\models T^{boole}$. If $\mathcal{M} \models 
 \beta'([\phi'_1(a,f)],\dots,[\phi'_{m'}(a,f)])$ for some tuple $a$ where $f$ is a tuple from $\mathcal{M}$, let $w_j:=[\phi'_j(a,f)]$. Then $T^{boole}\models \beta'(w_1,\dots,w_{m'})$. On the other hand, if $b_1,\dots,b_m$ form a partition of unity in $B$, $b_j\leq [\exists x \phi'_j(x,f)]$ for all $j\leq m$ where $f$ is a tuple from $\mathcal{M}$ and $\beta'(b_1,\dots,b_m)$ holds, then by Axiom 2.2(4)  
 $\mathcal{M}\models \exists x \beta([\phi_1(x,f)],\dots,[\phi_m(x,f)])$, and $(*)$ follows.
 
 Now let $\Phi$ be the set of formulas of the  form
  \[\beta( [\phi_1(x)],\ldots,[\phi_k(x)],u) )\]
  where $\phi_1(x),\ldots,\phi_k(x)$ are a partition of unity in $T$.
  It is easy to see that  $\Phi$ is  closed under Boolean combinations.     It follows immediately from the above argument that if $\phi(x,y) \in \Phi$
  where $y$ is a single variable, then $(\exists y) \phi$ is $T^{boole}$-equivalent to a formula in $\Phi$.  Quantifier-elimination for $T^{boole}$ follows.
  Stable embeddedness, and the statement on induced structure, follow from QE since there are no function symbols from $B$ into $\mathcal{M}$.
 QE includes the case of zero variables, so that sentences of $T^{boole}$  are   equivalent to sentences of $(B,[\psi])_{\psi}$ with $\psi$ varying over sentences of $L$.

 \eprf 
 
 Of course, if $T$ eliminates quantifiers, it suffices to include the formulas $[\phi]$ for quantifier-free $\phi$.

\ssec{Theories of Boolean algebras}

 Let $AB _{\forall}$  be the theory of Boolean algebras $B$  in the language of Boolean algebras $L_{BA}$ augmented by 
    predicates $AT_n, n=1,2,\cdots$;  $AT_1$ is intended to denote the atoms, and $AT_n$ the  sums of $n$ atoms, so this is merely an extension by definition.
 Let $ABA$ be the model completion. 
 It is axiomatized by $AB_{\forall}$ along with axioms asserting that for any element $u$ of $AT_n$
 there exist atoms $a_1,\ldots,a_n$ with $u=\sum_{i\leq n} a_i$; and that for any $u >0$ there exists $a \in AT_1$ with $a\leq u$.
 Let $ABA_n$ assert that there are exactly $n$ atoms, while $ABA_\infty$ asserts there are infinitely many. We will denote the theory of atomless Boolean algebras in the language of Boolean algebras by $BA_{atomless}^{boole}$.

  Let {{ABAI}}$_{\forall}$ be the expansion of $AB_{\forall}$ obtained by expanding the Boolean algebra $B$ by adding a  predicate 
  $J$ for an ideal of $B$ containing all atoms, as well as a  sort  $\bar{B}=B/J$, and a function symbol for the natural map $B \to \bar{B}$
 with kernel $J$.     Let ${{ABAI}} $ be 
  ${{ABAI}}_{\forall}$ along with an additional axiom asserting that $B \to {\bar{B}}$ is surjective, and ${\bar{B}}$ is atomless.  
 Let ${ABAI}_\infty$ assert, in addition, that there are infinitely many atoms of $B$. We define $ABAI_n = ABA_n +({\bar{B}}~ \mathrm{is ~ atomless}).$ We note that if $B/J$ is atomless, then $J$ contains all the atoms of $B$.
  
  \begin{lem}
  \label{mixed-basic-BA}   
  \begin{enumerate}
  \item $BA_{atomless}^{boole}$, ${{ABAI}}_\infty$ and each $ABAI_n$ is complete and eliminates quantifiers.   
  \item    
In a model of ${{ABAI}}$, $B/J$ is an atomless Boolean algebra, embedded and stably embedded in  $(B,J)$.
  \item  Let ${{ABAI}}^!$ be the class of models of ${{ABAI}}$ where $J$ is generated by the atoms of $B$.
Then any finite set of sentences consistent with ${{ABAI}}$ admits a model in $ABAI^!$.    
\item \label{mixed-basic-BA-4} Similarly, $ABA = Th(\{\Pp(n):  n=1,2,\cdots\})$.
  \end{enumerate}\end{lem}
  \prf  
  
 (1) These results are well-known. For a proof for the case of  $BA_{atomless}^{boole}$ see \cite[Theorem 6.22, p.85]{poizat-MT}. For a proof for the theories $BA_{atomless}^{boole}$, ${{ABAI}}_\infty$ and each $ABAI_n$ see \cite{DM-boole}.
 
 (2) can  be proved directly from (1) or by using the criterion of  Lemma 1 in the Appendix of \cite{ChH}:  an $\emptyset$-definable $D$ is embedded and stably embedded in $T$ if $T$ has a saturated model $M$ such that $Aut(M) \to Aut(D^M)$ is surjective. This also gives an alternative proof of (1).  
 
 Let $T={{ABAI}}_\infty$ or  $ABAI_n$.     Let $M,N$ be   saturated models of $T$ of the same cardinality.  Let $\bar{B}$ denote $B/J$, and 
 for $x \in B$ let $\bar{x}$ denote the class of $x$ modulo $J$.  
Let $G: \bar{B}(M) \to \bar{B}(N)$ be any Boolean algebra isomorphism ($G$ exists by completeness of the theory of atomless Boolean algebras.)
  Further let  $f: A \to B$ be an isomorphism of finite subalgebras of $M,N$ respectively, such that $G(\bar{a}) =\bar{f(a)}$ for $a \in A$.  
   Then we can easily extend $f$ to a bigger domain preserving these properties, 
  analyzing $1$-types over $A$, or by taking $M,N$ countable and using Stone duality.  Thus a back-and-forth argument provides
  an isomorphism $M \to N$ extending $G$ and $f$.  
  
  (3)  In other words, ${ABAI} = Th(ABAI^!)$.  This is clear from completeness of ${{ABAI}}_\infty$.
  
  (4) This follows from completeness of $ABA_\infty$.
  
  \eprf

 \ssec{Expansions of $T^{boole}$:  atomic, atomless and atomic-by-atomless}  \label{Tboole+}  
 
 We will consider three extensions of $T^{boole}$,  one in an expanded language.
 \begin{enumerate}
 \item   $T^{boole} _{at} := T^{boole} + ABA$, asserting in addition to $T^{boole}$  that the Boolean algebra $B$ is atomic.
 We will shorten the notation to $T_{at}$.

\item   $T^{boole}_{atomless} $ asserting  in addition to $T^{boole}$ that $B$ is atomless.   

 \item  \label{Lmixed} Let $L_{mixed,0} $ be $L^{boole}$ expanded by a predicate $J \subset B$.
 $T_{mixed,0}$ includes $T^{boole}$    and  asserts in addition that $B$ is atomic with infinitely many atoms, that  $J$ is a proper ideal of $B$
 with  is atomless. 
    
    Note that it follows that $J$ includes all atoms of $B$.   
    
  $L_{mixed}$ will be an inessential expansion of $L_{mixed,0}$.  
 $L_{mixed,0}$ includes a sort $S_B$ for each sort $S$ of $L$.    We define $L_{mixed}$ by 
  adjoining to $L_{mixed,0}$ a sort $\bar{B}$ along with a map $\pi_B: B \to \bar{B}$; and for each sort $S$ of $L$, a  sort $S_{\bar{B}}$ 
  and a map  $\pi_S: S_B \to S_{\bar{B}}$.  
%

 Let   $\sim_J$ be the equivalence relation defined by 
$x\sim_Jy$ iff $[x\neq y] \in J$. 

 $T_{mixed}$ extends $T_{mixed,0}$ by the assertions 
that the maps $\pi_B$ and $\pi_S$ are surjective with kernel given by  $\sim_J$.

  \end{enumerate}
  
We will put together all the $\pi_S$ and denote them by $\pi$.

We'll often write "$\mod \, J$" as a shorthand for "$\mod \, \sim_J$", and in particular $M/J$ for $M/\sim_J$.
 
\begin{rem}  If  ${L}$ has  a two-element sort $\two$ with elements $0,1$, the Boolean algebra $B$ of (1,2), and   the algebras $B$ and ${\bar{B}}$ of (3), do not require separate treatment; they can be identified with the Boolean interpretation of $\two$.  
 \end{rem}

 By \thmref{FV} and \lemref{mixed-basic-BA}, all three admit quantifier-elimination.  If $T$ is complete, then so is $T_{atomless}^{boole}$, while $T_{mixed}$ and $T_{at}$ become complete upon
 specifying that there are infinitely many atoms.

 For each $L$-formula $\phi$, the function $[\phi]:  M^n \to B$ induces a function  $ (M/J)^n \to B/J$, that we will  
 denote   $[\bar{\phi}]$.    The defining equation is:
   \[ [\bar{\phi}] (a\mod J) = [\phi](a) \mod J \]
   
 We will denote both quotient maps $M\to \bar{M}$ and $B\to \bar{B}$ by $\pi$; specifying the domain and range when needed, but
    the domains are in different sorts so there is no danger of confusion.   
   
  \begin{lem}  
  \label{mixed-basic}     Let $M \models T_{mixed,0}$, and let 
 \[\tM = (M,\bar{M},B,\bar{B}, \pi: B\to \bar{B},  \pi:M \to \bar{M},    [\phi]: M^n \to B, [\bar{\phi}]: \bar{M}^n \to \bar{B})\]
 (where $\phi$ ranges over all formulas of $L$) be   the corresponding model of $T_{mixed}$. 
 Then:  
 \begin{enumerate}
  \item   $\bar{M}=M/J$ is a model of $T_{atomless}^{boole}$.
 \item 
 Quantifier-elimination holds for $\tM$.  
 \item       $M/J$ is embedded and stably embedded in $\tM$. 
 \end{enumerate}
   \end{lem}

\prf   

It is clear that (1) holds, $\bar{M}$ is a model of $T_{atomless}^{boole}$.  Thus we already have:

\smallskip 

(2.1) \ \  Quantifier elimination for $(\bar{M},\bar{B}, [\bar{\phi}]: \bar{M}^n \to \bar{B})$.

\smallskip

We now prove (2). In view of (2.1), it remains to eliminate quantifiers over the sorts of $M$; this reduces to eliminating 
one existential quantifier over $M$.  Consider a quantifier-free formula $\theta(x, u)$ with $x$ ranging over $M$.   We need to find a quantifier-free formula equivalent to $(\exists x)\theta(x,u)$.   Equivalently, if we add constant symbols $d$ to the language matching $u$, we must
find quantifier-free formula equivalent to $(\exists x)\theta(x,d)$.   We thus supress the reference to $u$ in the notation.

The variable $x$ can only appear via $\pi(x)$ and various $[\phi(x)]$; i.e. $\theta(x)$ is equivalent to 
\[ \theta'([\phi_1(x)],\ldots, [\phi_n(x)], \pi(x) )\]
where $x$ does not occur in $\theta'$.  
  So we can rewrite $\theta$ as:
\[(\exists b_1 \in B)\ldots (\exists b_n \in B)(\exists \bar{a} \in M/J) \bigwedge_{i=1}^n [\phi_i(x)]=b_i  \wedge \pi(x)=\bar{a} \wedge 
\theta'(b_1,\ldots,b_n,\bar{a}) \]
Let $\psi(b_1,\ldots,b_n,\bar{a}) := (\exists a)(\pi(a)=\bar{a} \wedge \bigwedge_{i=1}^n [\phi_i(a)] = b_i$).
Then $(\exists x) \theta$ is equivalent to
\[(\exists b_1 \in B)\ldots (\exists b_n \in B)(\exists \bar{a} \in M/J) \psi(b_1,\ldots,b_n,\bar{a})  \wedge 
\theta'(b_1,\ldots,b_n,\bar{a}). \]

We will show that $\psi$ is equivalent to a quantifier-free formula; given this, 
 the quantifiers on $b_1,\ldots,b_n,\bar{a}$ can  be eliminated using (i).

Thus we must decide (on the basis of a quantifier-free formula),
given an  element $\bar{a}$ of $M/J$ and elements $b_1,\ldots,b_n \in B$, 
 whether there exists $x \in M$ with $\pi(x)=\bar{a}$ and   $[\phi_i(x)] = b_i$ for $i=1,\ldots,n$.    
Here we may assume $\phi_1,\ldots,\phi_n$ are a partition of unity in $T$.

A necessary condition is that $[\bar{\phi}_i(\bar{a})] = \pi(b_i)$ for each $i$.  Conversely assuming this condition holds, lifting $\bar{a}$ arbitrarily to an element $a'$
of $M$, we have $\pi(a')=\bar{a}$;  and for some element $e$ of $J$, $[\phi_i(a')] = b_i (\mod e)$ for each $i$.  Using Axiom \ref{axioms1} (\ref{lg1}), find 
$a''$ such that $[\phi_i(a'')] \meet e= b_i \meet e$.  Next glue $a',a''$ together, i.e. find $a$ with $[a=a''] = e$ and $[a=a'] = 1-e$.  Then 
$\pi(a)=\bar{a}$ and $[\phi_i(a)] = b_i$, as required.  

(3)  
From the nature of the quantifier-elimination (2)  it is clear that $M/J$ is stably embedded, and that the induced structure is the specified one.  

\eprf

 \begin{example}[Feferman-Vaught Products]  \label{products}

 Let   $I$ be a set.  Let $M_i$ be a model of ${T}$.    Let $M:=\prod_{i\in I} M_i$; if the language has several sorts, we 
 consider $M$ with the same sorts, defining $D(M):= \prod_{i\in I} D(M_i)$.     
  We use the language of \ref{boolean}, with the Boolean algebra $B$ interpreted to be $\Pp(I)$.   For a formula $\phi(x)$ of $L$ (where $x$ is a tuple of variables ranging over some sorts of $L$), 
 the  interpretation of $[\phi]$ in $M$ is \[[\phi](m) = \{i\in I: M_i \models \phi(m(i)) \}\]
 for $m\in M^n$. Note that $M \models T_{at}$.  
 
   Note that $B$ is interpretable in the original sorts of $M$, provided 
 ${T} \models (\exists x)(\exists y) (x \neq y)$  for variables $x,y$ of some sort $S$ of $L$, the language of $T$.    In this case 
 an element $a \in \Pp(I)$ can be coded as the equalizer of a pair of elements $c,d$ in $S$, i.e. $a=\{i\in I: c(i)=d(i)\}$.
If we assume further that $L$ has two constant symbols for elements of $S$, then $B$ is not only interpretable but definable.    In any case we will
view $B$ as a sort in the language of the product.  
   
The Boolean algebra of $\prod_{i\in I} M_i$, namely $\Pp(I)$, is atomic.   
 If $I_{at}$ is the ideal generated by the atoms in $B(M)$,   then the reduced
product $\prod_{i\in I} M_i /I_{at}$ has atomless Boolean algebra, and becomes a model of $T_{{atomless}^{boole}}$.   
Hence if $I$ is infinite, with $J:=I_{at}$, $M$ becomes a model of $T_{mixed}$.   

\end{example}

\begin{lem}\label{fmp} \begin{enumerate}
 \item \label{fmp1}  Let $T$ be a theory, and let $T^*$ be any completion of $T_{at}$.  Then any finite set $\Sigma$ of sentences of $T^*$ is 
realized in a model $\prod_{i \in I} M_i$ as in \exref{products}, with $I$ finite.   

\item  \label{fmp2}  Any finite set $\Sigma$ of sentences consistent with $T_{mixed}$ admits a model where $J$ is finite; in particular $J$ is generated by the atoms of $B$. 

\item  \label{perfectmodJ}  
   Let   $M = (\prod_{i\in I} M_i, \Pp(I), J) \models T_{mixed}$,  
   with $J$ the ideal of finite subsets of $I$.   
   Let $D'$ be a    definable set in $M^n/J$ (parameters allowed), and let $D$ be the pullback to $M^n$.  Then {\em either} $D'$ is finite (and in fact a finite union of one-element sets defined over the same parameters), {\em or}     there exists an infinite  $I_1 \subset I$ and an   injective, continuous map   $x:2^{I_1} \to D \subset M^n$, written $\eta \mapsto x_\eta$,  
   such that $x_{\eta} J x_{\eta'}$ iff $\eta E_{fin} \eta'$, where $E_{fin}  $ is the equivalence relation stating that $\eta-\eta'$ has finite support. 
   
    Moreover, one can arrange that 
 for any permutation $\si$ of $[n]$, if $x_{\eta} J x_{\eta'}^\si$ then $x_{\eta} J x_{\eta'}$.

\item   \label{atomlessacl}    Let $x$ be a variable on the main sort $S$ of $T_{atomless}^{boole}$; then the quantifier
   $(\exists^\infty x)$ is definable.   Moreover
    for $a \subset M^k$, $M \models T_{atomless}^{boole}$, $acl(a) \meet S = dcl(a) \meet S$.
\end{enumerate}

\end{lem}

 \noindent{\bf Remarks.} 
 
 (i)  \lemref{fmp} (\ref{fmp1})  is a result of \cite{FV}, answering a question of Mostowski.   
 
  (ii) Continuity of $x$ in (\ref{perfectmodJ}) refers to the  pointwise convergence topologies on $2^{I_1}$ and on $\prod_{i\in I} M_i$; in fact the construction will show that
   $x_{\eta}(i)$  depends only on finitely many specified values of $x(i)$.
    
    (iii) \lemref{fmp} (\ref{atomlessacl})
 holds  whether or not $(\exists^\infty)$ is definable in $T$.   

\prf (1)   By inspection of the axioms of $T_{at}$, they hold in any product structure, whether $I$ is finite or not.  Moreover $T^*$ is axiomatized
by $T_{at}$ along with $T$ relativized to the Boolean algebra $B$, enriched by constant symbols $[\psi]$ for each sentence $\psi$ of $L$.  
In particular $\Sigma$ follows from $T_{at}$ along with a certain consequence $\Sigma'$ of $T^*$  concerning $B$ and these constants.   
Only finitely many such constants $c_1 = [\psi_1], \ldots, c_k = [\psi_k]$ enter into $\Sigma'$; by taking the atoms of the Boolean algebra generated
by the $\psi_i$, we may assume the $\psi_j$ form a partition of unity.  
By the finite model property for the theory of atomic Boolean algebras,
we may find a finite $I$ and interpretations of $c_1,\ldots,c_k$ in the Boolean algebra $\Pp(I)$ such that $\Sigma'$ holds.  If $c_j \neq 0$, then $[\psi_j] \neq 0$ so there exists a model
$M_j \models T + \psi_j$.   For $i \in I$, let $M_i = M_j$ where $i \in c_j$.   Then it is clear that $\Sigma$ holds in $\prod_{i \in I} M_i$. 

(2)  The proof for $T_{mixed}$ is similar.  (It reduces to showing the same for  ${ABAI}_{\infty}$ itself, i.e. that $ {ABAI}_{\infty}$ is the theory of 
all pairs $(B,J)$ with $J$ a finite ideal and $B/J$ atomless.)

(3) By the quantifier-elimination \thmref{FV},  $d \in D$ is defined by 
\[ \Theta( [\phi_1(d)],\cdots,[\phi_K(d)] )  \]
where $\Theta$ is a formula of $(B,J)$ and $\phi_1,\ldots,\phi_K$ formulas of $L$, forming a partition of unity.  Replace each $\phi_k(x) $ by 
the two  formulas
\[ \phi_k^1(x):= \phi_k(x) \& (\exists^{\leq 1} y)\phi_k(y)\] 
and 
\[\phi_k^2(x) :=  \phi_k(x) \& \neg (\exists^{\leq 1} y)\phi_k(y). \]   
 Rearranging indices,
 we may further assume that for $k \leq K_1$,   $\phi_k$ has at most one solution in any model of $T$ and any assignment  of the parameters, while for $K_1 < k \leq K$, $\phi_k$ 
 has  
 at least two solutions in any model of $T$ and assignment. 
  So far we have not changed $D$, but only the representation.  

Now partitioning $\Theta$ and $D$ into finitely many pieces, we may further assume 
that $\Theta(u_1,\ldots,u_k)$ implies, for each $k \leq K$, either $u_k\in J$ or $u_k \notin J$.    

In case   $\Theta$ implies $u_k \in J$ for each $k >K_1$, 
 it is clear that $D/J$ is finite; in fact, having partitioned thus far, $D/J$ has at most one element.   

Assume therefore that on the contrary $\Theta \vdash u_{k_2} \notin J$ for some $k_2>K_1$.   Pick $a \in D$, and let $I_1 = [\phi_{k_2}(a)]$, an infinite subset of $I$.   
For $i \in I_1$ choose  $c_i \neq d_i \in  \phi_{k_2}(M_i)$.  For $\eta \in 2^{I_1}$, let $a_{\eta} (i) = a(i)$ if $i \notin [\phi_{k_2}(a)]$,   $a_{\eta}(i) = c(i)$ 
if $i \in [\phi_{k_2}(a)]$ and $\eta(i)=0$, and  $a_{\eta}(i) = d(i)$ 
if $i \in [\phi_{k_2}(a)]$ and $\eta(i)=1$.   Then $\eta \mapsto x_\eta=a_\eta$ is as promised.  Note that $[\phi_k(a_\eta)] = [\phi_k(a)]$ for each $k \leq K$.

For the `moreover' part of the statement we argue as follows. In case $\phi_{k_2}$ has $>n!$ solutions in each $M_i$, $i \in I_1$, it suffices, having picked $c_i$,  to pick   $d_i$ to differ from any of the  at most $n!$ permutations of $d_i$.  Since we only know that  $\phi_{k_2}$ has at least two solutions in each $M_i, i \in I_1$, we use a trick.
Let $l$ be an integer with $2^l > n!$, and divide $I_1$ into blocks of $l$ elements; so that $I_1 = [l] \times I_1'$.  We will  define $c_i,d_i$ not one index $i$ at a time but
in blocks of $l$ indices.   
For $i' \in I_1'$, let $M_{i'} = \prod_{j \in [l]} M_{(j,i')}$
be the product of the $l$ elements of the $i'$'th block.   Let $I' = (I \setminus I_1) \union I_1'$ and $i'=i$ for $i \notin I_1$.  
Let $\phi_{k_2}'$ be the product over the block of $\phi_{k_2}$.    Then $\prod_{i \in I} M_i = \prod_{i' \in I'} M_{i'}$.   But now  $\phi_{k_2}'$ has at least $2^l > n!$
solutions in each $M_{i'}$, so we can choose $c_{i'},d_{i'}$ so that no $Sym(n)$-permutation of $c_{i'}$ equals $d_{i'}$.   This concludes the proof of the `moreover' part.

(4)   
 Similar to the first part of (\ref{perfectmodJ}).     We can also view $S$ as a sort of $T_{mixed}$ (the main sort modulo $J$.)  We  have an $a$-definable set $D=D_a$ and we wish to find conditions on $a$ for finiteness of $D$.  
  We may again write $y \in D$ iff $\Theta( [\phi_1(y)],\cdots,[\phi_K(y)] ) $, where now the $\phi_i$ depend upon the  parameter $a$.
 The analysis of (\ref{perfectmodJ}) shows uniformly in $a$ that either $D_a$ is infinite,  or $D_a$ decomposes explicitly into
 a bounded finite number of $a$-definable sets with a single element.    Hence algebraic and definable closures coincide on $S$. 

\eprf

Recall that a theory $T$ admits 
  weak elimination of imaginaries iff  for any imaginary element $e$ in any model of $T$, we have
 $e \in dcl(e_1,\ldots,e_k)$ for some `real'  $e_1,\ldots,e_k \in acl(e)$. 
  Equivalently, any imaginary $e$ is equi-definable with a finite set
 of tuples real elements. For more details on stable embeddedness see e.g. \cite{CH} (2.3); for weak elimination of imaginaries see e.g. 
 \cite{casanovas-farre}.  Our main goal will be the proof of weak elimination of imaginaries for  $T_{at}$,  $T_{atomless}^{boole}$, and $T_{mixed}$, given elimination of imaginaries for $T$. This will be achieved in  Theorems \ref{1}, \ref{2}, and \ref{3}, respectively.  The remark below explains why assuming full EI for $T$ is necessary.

 \begin{rem}
 Weak elimination of imaginaries is preserved under finite products (see \lemref{I-finite}), but not under infinite products.

For instance, let $4$ denote a four-element pure set and ${4 \choose 2}$ the set of $2$-element subsets of $4$. Then $4 ^\Nn$ has ${4 \choose 2}^\Nn$ as an imaginary sort. Let $e \in {4 \choose 2}^\Nn$, viewed as an imaginary of $4^\Nn$. Then no element of  $4 ^\Nn$ is algebraic  over $e$, 
  so $4^\Nn$ does not admit weak elimination of imaginaries.  Nor does the theory with $\clubsuit$
 obtained by adding an ideal element, as in \remref{clubrem}.     (We leave these statements to the reader, noting that
  if $G$ acts on $X$ without fixed points, 
 then $G^\Nn$ acts on $X^\Nn$ without  finite orbits.   Apply this to $G=Aut(4;\{1,2\})$ acting on $4$.)

 \end{rem}

\section{Imaginaries in infinite products} \label{sec:atomic}

 \begin{thm}  \label{1}    Assume ${T}$ admits EI, and each sort of $T$ has a $0$-definable element.  Then $T_{at}$  admits weak EI.   \end{thm}
Following some preliminaries, we will  prove the theorem for the very special case of the theory $T=SET_n$ of an $n$-element set,   with a single distinguished point (\lemref{setn}.)  
Here we will use liaison groups in a somewhat new setting, and make use of results of descriptive set theory.   
 We will then prove \thmref{1} in general, by an argument reducing to enrichments of $(SET_n)_{at}$.

We begin with the case of finite index set $I$.

\begin{lem} \label{I-finite}  [$\clubsuit$]   Assume $I$ is finite, and $M_i$ has weak EI for $i \in I$.  Consider   $M=\prod_{i\in I} M_i$ as a structure as in \ref{products}.
Then $M$ admits weak EI.   \end{lem}
  \prf Note that the structure of Example \ref{products} includes, in this case, the finite Boolean algebra $\Pp(I)$.  We have in general:
  
  \claim{}  Consider an expansion $(M,c)$ of a structure $M$  by naming a tuple $c$ of  algebraic (real) elements.  If $(M,c)$ has weak EI, then so does $M$.
  
  \prf   For any imaginary $e$ of $M$ we have,
  using weak EI in $(M,c)$, that 
  $e \in dcl(e_1,\ldots,e_k,c)$ for some `real'  $e_1,\ldots,e_k \in acl(e)$; since $c$ is also a real element in $acl(e)$, the condition for weak EI in $M$ is satisfied. \eprf
  
   Thus returning to the lemma, we may as well add names to the elements of $I$ (identified with the atoms of $B$.)  
 In this formulation the lemma reduces to the case $|I|=2$, $M=M_1 \times M_2$.   We will also view a finite product of real sorts as a real sort itself.
 
 Let $Z$ be a definable subset of $M^n$ , with canonical parameter $e \in M^{eq}$; 
 we will find a finite subset $W_i$ of $M_i$ such that $e$ is interdefinable with a code for a certain  subset $W$ of $W_1 \times W_2$.   
 We note first that $Z$ is {\em canonically} a finite union of some rectangles $X_i \times Y_j$ with $X_1,\cdots,X_k$   definable in $M_1$
 and $Y_1,\ldots,Y_l$ definable in $M_2$.    Namely, let $\mathcal{Y}_2 = \{Z(a):  a \in M_1\}$; where 
 $Z(a) = \{y: (a,y) \in Z\}$.   Then $\mathcal{Y}_2$ is a {\em finite} set of $M_2$-definable subsets of $M_2$.  
 (To see this, it suffices to saturate $(M_1,M_2, M_1 \times M_2)$ and then note that $\mathcal{Y}_2$ is bounded.   This is because an automorphism of $M_1$ extends to an automorphism of
$M_1 \times M_2$ that fixes $M_2$.  So $Z(a)=Z(b)$ if $tp_{M_1}(a)=tp_{M_1}(b)$).

 Define $\mathcal{Y}_1$ dually.
 Then $Z$ is a union of some of the `blocks' $U_1 \times U_2$, with $U_i \in \mathcal{Y}_i$.   This property continues to hold if $\mathcal{Y}_i$ 
 is replaced by a larger finite set.  Using
  weak EI in $M_i$,  each element of $\mathcal{Y}_i$ is interdefinable with some finite subset of some (`real') sort of $M_i$; hence 
 we may assume $\mathcal{Y}_i = \{U_b:  b \in W_i\}$, where $W_i$ is a finite subset of some sort in $M_i$, and $(U_w)_w$ is a 0-definable family
 of definable sets in $M_i$.  Let $W= \{(w_1,w_2): w_i \in W_i, U_{w_1} \times U_{w_2} \subset Z \}$.    Then  $W$ is definable from $Z$, and  
 \[ Z = \bigcup \{ U_{w_1} \times U_{w_2} :  (w_1,w_2) \in W \}.\]  
 
 If $S_i$ is the sort of the $U_{w_i}$'s, we can write 
 \[U_{w_1} \times U_{w_2}  = (U_{w_1} \times S_2) \meet (S_1 \times U_{w_2}) \]
 (and $U_{w_i}$ is the projection of $U_{w_1} \times U_{w_2}$ to the $i$'th coordinate).
 So it remains to show that sets such as $U_w \times S_2$ have a real code algebraic over their canonical code $e_1$; in $M_1$, such a real code $c_1$ does exist; we now use $\clubsuit$ to find a $0$-definable $c_2$ of $M_2$; then $(c_1,c_2)$ is a real code for $U_2 \times S_2$, algebraic over $e_1$.
 
 \eprf

 \begin{rem} \label{I-finite2}   Assume $M_1,\ldots,M_n$ each  admit EI.   For $i<n$ assume also that $M_i$ has a definable linear ordering or more generally rigid finite sets, i.e. no nontrivial automorphism
 of an elementary extension of $M_i$ has a finite orbit.  Then $\prod_{i \leq n} M_i$ admits EI.  (This follows from \lemref{I-finite}.)  \end{rem}
   
  \begin{lem} \label{u0}  
  \begin{enumerate} \item Let $I$ be an infinite set.  Then $\Pp(I)$ is  an $\aleph_0$-saturated model of the theory of atomic Boolean algebras. 
  \item   Assume $T \models (\exists ^{\leq n} x)(x=x)$.  
   Then in any model $(N,B) \models T_{at}$  there exist finitely many parameters $c$ in $N$ such that $N \subset \dcl(B,c)$. 
  \item    Assume $T \models (\exists ^{\leq n} x)(x=x)$.    Let $M \models T$ and form the power $N =M^I = \prod_{i \in I} M$ 
   as in \exref{products}, with Boolean algebra $B=\Pp(I)$.     Then $N$ 
    is an $\aleph_0$-saturated model of $T_{at}$.
  \end{enumerate}
  \end{lem}
  
  \prf  (1)  Let $J$ be the ideal of $B=\Pp(I)$ generated by the atoms.  Then $\bar{B}=B/J$ is atomless and infinite.  Any finite subset of $B$ is contained
  in a finitely generated subalgebra $S$ with disjoint generators $b_1,\cdots,b_k, a_1,\cdots,a_l$, where the $a_i$ are atoms and no $b_j$ is a finite union of atoms in $B$.
  Note that modulo $J$, the elements $b_1,\ldots,b_k$ are precisely the atoms of $\bar{B}(S)$.  A $1$-type $tp(c/b_1,\cdots,b_k, a_1,\cdots,a_l)$ is determined by the set of $a_i$ below
 $c$, and whether $c \meet b_j$ is a finite union of atoms, or the complement of that, or neither.  It is clear that all such 1-types are realized.
  
  (2)   
  Inductively find $a_1,\ldots,a_n$ such that for each $j \leq n$, $[\bigwedge_{i<j} a_i \neq a_j]= [(\exists x)(\bigvee_{i<j} a_i \neq x]$.   Then
  $[(\exists x)(\bigvee_{i \leq n} a_i \neq x]=0$.   So any $c  \in N$ is definable from $a_1,\ldots,a_n$ and  the $n$ elements $[c=a_i]$ of $B$.

  (3) There are finitely many   formulas $\phi(x_1,\ldots,x_n)$ of $Th(M)$, so by   \thmref{FV}, for any $a_1,\ldots,a_n \in M^I$,
  the  induced structure   on $B$ from $(M^I,a_1,\ldots,a_n)$ is just $(B,e)_{e \in E}$ for some finite set $E$.    Hence $\aleph_0$-saturation of $M^I$ follows from (1) and (2).

  \eprf

  We will also use  Wencel's theorem from  \cite{wencel} that we quote as a Proposition:    
   \begin{prop} \cite{wencel} \label{wencel} The theory ABA of atomic Boolean algebras  admits weak EI.
  \end{prop} 
  \begin{cor}\label{smooth} Let $B =2^{\Nn}$, an atomic Boolean algebra,  and let $E$ be a definable equivalence relation on  a definable $D \subset B^n$.   Then $E$ is smooth.  
  \end{cor}
  \prf  
  If $D$ is any definable set of $B^n$, and $E$ any definable equivalence relation, then $D,E$ are Borel and $E$ is a smooth equivalence
  relation on $D$.  The Borelness of $D,E$ is immediate from quantifier-elimination, while the smoothness follows from    weak elimination of imaginaries in $B$:  it provides a definable (hence Borel) map $q$
  from $D$ to the  space of finite subsets of $B^m$ (for appropriate $m$), such that $q(d)=q(d')$ iff $dEd'$.   
  \eprf

  We will see how to propagate this smoothness from the Boolean algebra to product structures.     
  In \cite{bh}, a completely general direct proof is given that definable equivalence relations on countable products of countable structures are
  smooth; giving moreover an alternative route to a proof of \thmref{1} that avoids special analysis of the Boolean algebra.   However the model-theoretic methods used here with internality and liaison groups appear to us to have independent interest.

 \ssec{Definable groups and automorphism groups}  \label{definablegroupsinB} Consider the Boolean algebra $B=\Pp(I)$, where $I$ is a countably infinite set.  For any finite structure
 $F$, $F^I$ is interpretable in $B$:  take first the case $|F|=2^m$, and let $\{c_1,\ldots,c_{2^m}\}$ be an enumeration of $F$.  We have identifications
 $B^m = (2^I)^m = (2^m)^I = F^I$.   The general case also follows by allowing $\{c_1,\ldots,c_{2^m}\}$ to be a listing with repetitions, or by viewing $F$ as a
 definable subset in some richer structure of size $2^m$.  
 
     In particular, for a finite group ${G}$, the group  ${G}^I$ is interpretable in $B$.  To see this in another way,
  one can think of $B$ as the commutative ring $\Ff_2^I$, embed a group ${G}$ in $GL_n(\Ff_2)$ and hence embed ${G}^I$ in $GL_n(\Ff_2^I) = GL_n(\Ff_2)^I$.

 \begin{question}   Is every interpretable group over $B=2^I$ definably isogenous to a subgroup $GL_n(\Ff_2^I)$?   More precisely, let $G$ be a group interpretable in $B$.  Does there exists a definable finite index subgroup $G_1$ and a definable homomorphism $h: G_1 \to GL_n(\Ff_2^I)$, for some $n$, with finite kernel?\end{question}

 $G^I$ has the structure of  a separable compact group.  
 
 It follows   from quantifier-elimination for $B$,  that all definable subsets of $B^m$ (with parameters) are Borel.  
 In general, definable subsets of ${G}^I$ need not be closed.  (Witness the elements supported on one atom.)  
 
 To bridge the gap between Borel and closed, we will use a result of Becker and Kechris.  
  Taking into account  the Harrington-Kechris-Louveau theorem \cite{hkl} or \cite{bk} 3.4.3, Theorem 3.4.5 of \cite{bk} asserts:    
\noindent{\em  Let $\bX$ be a Polish space, $\G$ a group acting on $\bX$ by homeomorphisms, and let $R$ be the equivalence relation on $\bX$ whose classes are the $\G$-orbits.     If there is a dense orbit and $R \subset \bX^2$ is meager,  then $R$ is not smooth.}
 
In a significant special case this can be sharpened:  
 \begin{cor} \label{meager} Let $H$ be Borel subgroup of a   Polish group $G$, with closure $\bH$.   Assume  the  coset equivalence $xH=yH$ on   $G$ is smooth. 
Then $H=\bH$.   
 \end{cor}

\prf  
 
Consider $\G=H$ acting on $\bH$ on the left.    $H$ itself is a dense orbit of this action.
Since $Hx=Hy$  is smooth on $G$, the restriction to $\bX$ is still smooth.        
The orbits of $\G$ are $H$-classes, and the orbit equivalence relation is   \[R=\{(x,y) \in \bH^2:  y x \inv \in H\}.\]
If $H$ is meager in $\bH$, it is easy to see that $R$ is meager in $\bH^2$; but  this contradicts \cite[Theorem 3.4.5]{bk}.  Hence $H$ is not meager.  
    By \cite[Theorem 9.9]{kechris}, 
  $H=H H \inv$ contains an open neighborhood $U$ of the identity of $\bH$ ; hence $H=HU = \bigcup_{h \in H} hU$ is open in $\bH$, and thus also closed in $\bH$, so $H=\bH$. 
   \eprf

 \begin{lem}  \label{dst} Let $H$ be a subgroup of $G^I$, definable with parameters in $\Pp(I)$.    Then $H$ is closed.    \end{lem}
  \prf We may assume $I$ is countable by a L\"owenheim-Skolem argument.   (If $\fm$ is a countable elementary submodel of the relevant structure, and $H$ is not closed in $G^I$,  
  then for some $n \in G^I \setminus H$, for any finite $I_0 \subset I$, there exists $h \in H$ with $I_0 \subset [h=n]$; the same remains true in the
  elementary submodel $\fm$: note that $(G^I)(\fm) \subset G(\fm)^{I(\fm)}$ naturally; so 
    $H(\fm)$ will not be closed in $G^{I(\fm)}$ either.) 

  From definability and quantifier-elimination  we see that  $H$ is Borel.   By \corref{smooth} the $H$-coset equivalence relation is smooth.  Let $\bH$ be the closure of $H$.
  From \corref{meager} we conclude that $H=\bH$.
  \eprf
  
     We will  need a extension of \lemref{dst}.   
       Let $D=\coprod_{i \in I} \Om_i$, with $|\Om_i| = m  \in \Nn$ for all $i\in I$; and let $N_0$  be the structure $N_0=
    (D,I,f)$,  with sorts $D,I$ and with $f(x)=i$ for $x \in \Om_i$.     Then $Aut(N_0)$ has $Sym(I)$ as a quotient group.
    Let $G=Aut(N_0)$, $H = Aut(N_0/I) = \prod_{i \in I} Sym(\Om_i)$.    Note $G=Aut(N_0)$ is naturally a Polish group.

\begin{lem}  \label{closedgroup-r}  Let  $\widehat{G}$ be a Borel subgroup of $G$. Assume $\widehat{G}$ projects onto $Sym(I) = G/H$, with kernel $\widehat{H}$.
\begin{enumerate}   
 \item     Assume $\widehat{H}$   is  definable with parameters in the structure $\prod_{i \in I} \Om_i$.    Then $\widehat{G}$ is a closed subgroup of $G$.
  \item  Assume   $\widehat{G}$  is a closed subgroup of $G$.  Then $\widehat{G}$ 
 has a closed finite index subgroup of the form $Aut(N_0/e^*)$ for an appropriate function $e^*$, such that $e^*(i) \in  \Om_i^{eq} $.  
 \end{enumerate}
  \end{lem}   
  
  \prf   
    (1)  By \secref{definablegroupsinB},  $H$ is definable in $B$.  By assumption, $\widehat{H}$ is a definable subgroup of $H$.  
     By \corref{smooth},  the $\widehat{H}$-coset relation on $H$ is smooth, i.e.
    there exists a Borel map $f$ on $H$ into a Polish space, such that $f(x)=f(y)$ iff $x\widehat{H}=y\widehat{H}$.
  By assumption, $\widehat{G} H = G$, so for any $x \in G$ there exists $\hat{x} \in H$ with $\hat{x} \widehat{G} = x \widehat{G}$.
   Then the function $\hat{f}(x) :=   f(\hat{x})$ is well-defined, analytic and hence, Borel, and maps $G$ to a Polish space, with kernel $\widehat{G}$.   In other words  the  $\widehat{G}$-coset equivalence relation is smooth.    
 By Corollary \ref{meager}  $\widehat{G}$ is closed.    
 
(2)    Let $N$ be the expansion of $(D,I,f)$ by all relations invariant 
  under $\widehat{G}$.   Since $\widehat{G}$ is closed, we have $\widehat{G} =Aut(N)$.  
  Note that $N$ is $\aleph_0$-categorical:  any $k$-tuple is contained in the union of at most $k$ of the sets $\Om_i$, so 
  the number of orbits on $k$-tuples is at most $(m k)^k$).   
  
  Within $N$,   $I$ is a trivial strongly minimal set.  Strong minimality and triviality for the induced structure on $I$ are clear, since every definable relation is $Sym(I)$- invariant and hence 
definable in pure equality.  In general, to say that a definable set $I$ is strongly minimal within some structure $N$ is stronger than 
just saying this for the induced structure, since you need every $N$-definable set to be finite or cofinite.  
But in the present case $N$ is contained in $acl(I)$,  so any type of elements of $N$ over $I$ is isolated and hence $I$ is stably embedded in $N$.
Stable embeddedness can also be deduced abstractly from the surjectivity of $Aut(N) \to Aut(I$).
  
So $N$ is an almost strongly
  minimal structure (contained in the algebraic closure of a strongly minimal definable set), and $N$ has Morley rank $1$.     Pick $c \in I$, and let $d$ enumerate $f \inv(c)$.      
  In $N^{eq}$, define the set  $A := \dcl(d) \meet \acl(0)$.    Then $tp(d/A)$ is stationary.   Clearly $A$ does not depend on the specific enumeration $d$.  It does not depend on $c$ either, since $tp(c/0)$ is stationary, so $tp(c/0)$ implies a complete type over $\acl(0) \supseteq A$.  Thus $A$ is an
  $Aut(N)$-invariant subset of $acl(0)$.  Since also $A \subset \dcl(d)$, by  either $\aleph_0$-categoricity (using finitely many types $tp(d,d')$) or
  $\omega$-stability (using the Morley rank and degree of $tp(d)$), there can be no strictly increasing chain of definable subsets of $A$, so   $A$
  is essentially finite, i.e. $A=\dcl(A_0)$ for some finite $A_0 \subset A$. 
  
    Moreover if $c_1,\ldots,c_k$ are distinct elements of $I$
  and $d_i$ an enumeration of $f \inv(c_i)$, then $d_1,\ldots,d_k$ are independent over $A$; so $tp(d_1,\ldots,d_k/A)$ is uniquely determined given $tp(d_1/A),\ldots,tp(d_k/A)$.  Hence 
  letting  $H_c = Aut_N(f \inv(c) / c,A)$, we see that    $\prod_{c \in I} H_c$ acts on $D$ by automorphisms over $I,A$.  The 
  homomorphism $\prod_{c \in I} H_c \to Aut(N/ A,I)$ is clearly surjective, with trivial kernel.   Thus
   $\prod_{c \in I} H_c = Aut(N/ A,I)$.   Now we have a homomorphism $Aut(N/I) \to Sym(A)$, by restriction.  As above, the kernel $Aut(N/ A,I)$ has finite index since   $A$ is essentially finite.  Thus $Aut(N/ A,I)$   is a finite index subgroup of $\wH = Aut(N/I)$.    
   
The $H_c$-orbit  of $d$ is definable in $N$ over $Ac$; so there exists an $A$-definable function $e^*$ such that $e^*(c)$ codes the $H_c$-orbit
of $d$.  Now $e^*$ is $Aut(N/A)$ - invariant, so the finite-index subgroup $Aut(N/A)$ of $\widehat{G} =Aut(N)$ fixes the function $e^*$, so 
$Aut(N/A) \leq Aut(N_0/e^*)$. On the other hand $Aut(N_0/e^*,I) = \prod_{i \in I} H_i = Aut(N/A,I)$;  since the two groups map onto $Sym(I)$
and have the same kernel, the inclusion $Aut(N/A) \leq Aut(N_0/e^*)$ must be an equality $Aut(N/A) = Aut(N_0/e^*)$.

  \eprf

\begin{rem} 
  
  The  analogue in the atomless case would be the Boolean power of $\Om$ with respect to an atomless Boolean algebra; in this 
case  the conclusion of \lemref{closedgroup-r} can probably
 be strengthened from commensurability to equality, i.e. $\widehat{H}$ itself has the right form.  We will in any case not need this.  \end{rem}

Let $SET_{\leq n} $ be the theory, in the language with one constant symbol, of a nonempty set $\Om$ with at most $n$ elements;
augmented with imaginary sorts coding subsets of $\Om^k$ for $k \leq n$.
 Thus the models of $SET_{\leq n} $ are finite altogether.  The constant symbol serves to meet the condition $\clubsuit$
 (note that for any $k$ there is a 0-definable subset of $\Om^k$.){

    It is easy to see that $SET_{\leq n} $ eliminates imaginaries to this finite language (even without the constant symbol).    Let $SET_n$ be the completion where the home sort has exactly $n$ elements.

 \begin{lem}\label{setn}   $(SET_{\leq n})_{at}$  admits weak EI.   \end{lem}
 
\prf  A model  $M$ of $(SET_{\leq n})_{at}$ can be identified with a product $\prod_{1 \leq i \leq n} M_i$, with $M_i \models (SET_i)_{at}$.
Hence using \lemref{I-finite}, it suffices to prove weak EI for  $T=SET_{n}$ instead.

By \lemref{u0} (3) $T_{at}$ has a countably saturated model $M=\Om^I$ with $|\Om|=n$, $|I|=\aleph_0$.   We have $B=B(M) = \Pp(I)$.   
  We observed in  \lemref{u0} (2)  that $M$ is internal to $B$,
   i.e. $M \subset dcl(B,c_1,\ldots,c_n)$ for appropriate elements $c_1,\ldots,c_n$.   
   Recalling also the stable embeddedness of $B$   (\thmref{FV}), 
   we   have a definable liaison group    $H=Aut(M/B)$ 
   (cf. \cite[Appendix B]{bedlewo}). 
    Let $G=Aut(M) \cong Sym(I) \sd H$.
   The same $G$ is also the automorphism group of the countable structure $\coprod_{i \in I} \Om_i$, where $|\Om_i|=|\Om|$.   We use 
   the latter to endow $G$ with a topology of pointwise convergence, and in particular a Borel structure.
    
   Let $e \in M^{eq}$ be an imaginary.   Let $A=dcl(e) \meet B^{eq}$.    We have $e \in dcl(d_1,\ldots,d_m)$ for some $d_1,\ldots,d_m \in M$.
   By \thmref{FV}, the induced structure on $B$  from $M(d_1,\ldots,d_m)$  is the Boolean algebra structure along with constants
   $[\psi(d_1,\ldots,d_m)]$ where $\psi$ is an $m$-ary formula of $T$.  Since there are only finitely many such formulas, the induced structure on $B$
   from $(M,e)$ 
   is a reduct of $B$ enriched with finitely many constants; in other words it is the Boolean algebra enriched with an imaginary constant $A$.  
   By weak elimination of imaginaries in $B$ (\propref{wencel}),  
   there exist elements $b_1,\ldots,b_k \in B$
   (that we may take to form a Boolean partition of $1$)  such that each $b_i \in acl(A)$, and $A \subset dcl(b_1,\ldots,b_k)$.
     By stable embeddedness of $B$,
   $tp(e/A)$ implies $tp(e/B)$.     
   
    Now $B$ splits as the product of $k$ Boolean algebras $B_i$ (with $1_{B_i}=b_i$.)
   Likewise $M$ can be viewed as $M=\prod_{i=1}^k M_i$ where $M_i$ is the part corresponding to $B_i$.  
   Using wEI for finite products (\lemref{I-finite}), $\acl(e)$ includes imaginaries $e_i$ of $M_i$, with $e \in \dcl(e_1,\cdots,e_k)$.
   So we can reduce to the case $k=1$ and $b_1=1$.   In other words $A \subset \dcl(b_1)=\dcl(\emptyset)$; so $\dcl(e) \meet B^{eq} = dcl(0)_B$, or again, by stable embeddedness,
   $\widehat{G}:=Aut(M/e) $ maps onto $Aut(B)$.    Note also that  $\widehat{G}$ is a Borel subgroup of $G$; if $e=a/E$ for a definable equivalence
   relation $E$ on $M$, we have $g \in \widehat{G}$ iff $g(a) E a$.
   
    Let $\widehat{H}=Aut(M/B,e)$.  Since by \lemref{u0}(2) $M$ is $B$-internal, again by the theory of internality (\cite[Appendix B]{bedlewo}),  the liaison group $\widehat{H}$ is an $M$-definable (with parameters) subgroup of $G$, or an intersection of a descending chain of $M$-definable subgroups $\widehat{H}_j$. 
    Since $\widehat{H}$ is also definable over finitely many parameters $e$, by $\aleph_0$-categoricity of $M$ it is definable;
 $\widehat{H}= \widehat{H}_j$ for some $j$.  Viewing $M$ and $G$ as interpretable with parameters over $B$, we see  that $\widehat{H}$ is definable with parameters in $B$.    
    By \lemref{closedgroup-r},  $Aut(M/e)=\widehat{G}$    has a finite index closed subgroup of the form $Aut(M/e^*)$
    for some function $e^*$, such that $e^*(i) \in \Om_i^{eq}$.   Seen as an  element of $M$, we have  $e^*  \in \acl(e)$ (by homogeneity of $M$ and the finiteness of the orbit of $e^*$  under $\widehat{G}$). 
    Conversely 
 $Aut(M/e^*)$ fixes $e$, so $e \in \dcl(e^*)$.  This 
proves weak elimination of imaginaries for $(SET_n)_{at}$. 
 \eprf

Before proceeding, let us show Lemma \ref{setn} immediately yields the case of finite sets with arbitrary structure. 
This will not be used in the sequel since a somewhat more general version  will be required, namely enrichments of the infinite product that do not add structure to the Boolean algebra.

\begin{cor}\label{2.2} ($\clubsuit$).
If ${T} \models (\exists ^{\leq n} x)(x=x)$,   then   $T_{at}$  admits weak EI.
\end{cor}
\prf   
Let $T'$ be the reduct of ${T}$ to the language with one constant symbol, using the 0-definable element given by $\clubsuit$ to interpret that constant symbol.

Let $T''$ be an expansion of $T$ to the language with
$n$ new constant symbols $c_1,\ldots,c_n$, and let ${T}''$ assert that they   
exhaust the universe.    

Let $M$ be any model of $T_{at}$.  By
\lemref{u0} (2), $M$ can be expanded to a model $M''$ of $T''_{at}$.  Let $M'$ be the reduct of $M$
 to a model of $T'_{at}$.   Then $M''$  is an expansion of $M'$ by   $n$ constants.
 So the parametrically defined sets of $M',M''$ and hence also of $M$ coincide.  By  \lemref{setn}, each such set has a   canonical
  parameter consisting  of codes for finite subsets of tuples from $M'$; 
     this is still a   canonical parameter from the point of view of $M$,   proving weak EI for $M$.    \eprf

  We now move towards the proof of weak elimination of imaginaries for  the atomic case (cf. \thmref{1}). The tools that we develop in the following section will play an important role in the mixed case, Theorem \ref{3}, allowing us  to code a subset $Z$ contained in a single $J$-class.   
  
  \ssec{Definition of $E^Z$}  \label{ez}  Let $(M,B) \models T_{at}$.    Let $D$ be a sort of $T$, and let $Z=Z_c \subset D$ be an $M$-definable set in the language of $T_{at}$, with imaginary canonical parameter $c$.  
  Let $C = \acl(c) \meet M$, i.e. $C$ is the set of elements in `real' sorts of $(M,B)$ that are algebraic over the imaginary parameter $c$. 
  We will write $M_C$ for the expansion of $(M,B)$ obtained by naming the elements of $C$.

    Our plan is to  define   equivalence relations $E_a^Z$  of $T$ on $D$, uniformly for each atom $a$ of $B$.  We will show that for some $n$, we have $|D/E_a^Z| \leq n$ for each atom $a$ of $B$.   We will then 
let $E^Z$ be the equivalence relation on $D(M)$, with   $x E^Z y$ if and only if for all atoms $a$ of $B$, $xE_a y$;  and show that  
  $Z$ is a union of $E^Z$-classes.    This will essentially reduce the problem to the 
  case $D$    is `definably pro-finite',   i.e. the Boolean value of 
   $|D| \leq n$ is $1$, for some $n$.

  Let $a \in B, a'=1-a$.  We can identify $B$ with a product algebra $B_a \times B_{a'}$, with $a$ identified with $(1,0)$ and $a'$ with $(0,1)$.
  (So $B_a = \{y \in B: y \leq a\}$, with $1_{B_a}=a$  and $\neg y = a \m y$ for $y \in B_a$).   
  Then $M$ also decomposes definably as $M = M_a \times M_{a'}$, where $M_a = M / \Theta_a$, $x \Theta_a y$ iff $a \leq [x=y] $.  In case 
  $M=\prod_{i\in I} M_i$,   the factor $M_a$ can be identified with   $\prod_{i \in a} M_i$.  Note $M_a,M_{a'}$ are interpretable in $(M,a)$, and stably embedded and orthogonal to each other.
  (Note that the automorphism group of $(M,a)$ is just the product $Aut(M_a) \times Aut(M_{a'})$).  Thus also $D(M)  = D(M_a) \times D(M_{a'})$.
  For $d \in D(M_a)$,  let $Z(d) =\{d' \in D(M_{a'}):  (d,d') \in Z\}$.   
 
      We obtain from $Z$ alone a definable map $D(M_a) \to  W(M_{a'})$ for a certain imaginary sort $W$ 
   mapping $d \in D(M_a)$ to the code of $Z(d)$.   
  Since $M_a,M_{a'}$ are  orthogonal,  this map has finite image.  
 \footnote{In general, a definable map from one interpretable set to an orthogonal one must have finite image; this follows from the fact that the graph of the map is a finite union of rectangles.}
   Thus we canonically obtain from $Z$ finitely many pairwise disjoint subsets $Z_{a,1},\ldots,Z_{a,r}$ (for $M_a$) and imaginary elements
   $e_1,\ldots,e_r$ (for $M_{a'}$) such that $Z_{a,i} \mapsto e_i$; but the orderings of the $Z_{a,i}$ and the $e_i$ are not canonical.

   Let $E_a^Z$ be the equivalence relation on $D(M_a)$ defined by  $d E_a^Z d'$ iff $Z(d) = Z(d')$.   By the above discussion,
   $E_a^Z$ has only finitely many equivalence classes.  Since this is true in any model of $T_{at}$, by compactness, the number is bounded by some $n$.

   By elimination of imaginaries in $T$, $E_a^Z$ is coded by some tuple $\g_a$ from $M_a$, definable uniformly in an atom $a$ of $B$.   By Axioms 2.2(4)  of $T_{boole}$,
   there exists an element $\g$ with $\g =\g_a$ (modulo $a$) for each $a$.   Alternatively, this is clearly true in $\prod_{i \in I}M_i$ with finite $I$, hence by \lemref{fmp}(\ref{fmp1}) 
   it holds in $M$.    
       We defined $\g$ canonically from $Z$, hence by the definition of $c$ and $C$ we have $\g \in C$.   So $E_a^Z$ is definable uniformly in $a$ over $C$.
 
   Define $x E^Z y$ if and only if for all atoms $a$ of $B$, $xE_a^Z y$.    So $E^Z$ is $C$-definable.
   
   \begin{rem} \label{ez2}  The construction of $E^Z$ makes sense for $Z$ defined in any expansion $T_{at+}$ of $T_{at}$, provided the orthogonality
   of $M_a$ with $M_{a'}$ mentioned above continues to hold for all atoms $a$ of $B$.     \end{rem}
   
For now we return to the case where $Z$ is defined in an expansion by constants of $T_{at}$.
   \begin{lem}  \label{oneatomatatime}     
   $Z$ is a union of $E^Z$-classes.
\end{lem}

\prf  This is an elementary statement, so by \lemref{fmp}(\ref{fmp1}) it suffices to prove it in product structures $\prod_{i \in I} M_i$ with $M_i \models T$ and $I$ finite.  In such
models the statement is clear:  if $x E y$, i.e. $xE_i y$ for each $i \in I$,   one can move from $x$ to $y$ changing one coordinate at a time to see that $x \in Z$ iff $y \in Z$.
\eprf

 \prf[Proof of \thmref{1}]    
 Let $M$ be a sufficiently saturated and homogeneous  model of  $T_{at}$.   
  Let $Z$ be an $M$-definable subset of $M^n$.  
 We need to find a canonical code for $Z$.     Recall the definition of $C=\acl(\ulcorner Z \urcorner) \meet M$ and $M_C$ from \secref{ez}, as well as the equivalence relations $E_a$ and $E$.
 
 By EI in $T$, we may view $D/E_a$ as a finite subset $F_a$ of some sort $S$ of $T$.  We have $|F_a| < n$ for some $n$.  Using $\clubsuit$,
 we can throw in a $0$-definable (in $T$) element $\om \in S$   into $F_a$.  Then $D/E$ embeds definably into the sort $S^M$ (with its Boolean $T_{at}$-interpreation). Let $Q=(D/E,B)$, with the structure of a  model of $(SET_{\leq n})_{at}$:   the Boolean algebra structure on $B$, the maps
 $[x=y]:  (D/E)^2 \to B$ and $[x=\om]: (D/E) \to B$.  
 So $Q$ is interpretable in $M_C$.  
  Of course, the induced structure on $Q$  from $M_C$ may be richer; 
 but we have:
 
 \claim{}   Let $W \subset Q^k$ be $M$-definable.  Then $W$ is definable with parameters in $Q$.
 
 \prf   By \lemref{u0} (2), $D/E^Z$ is internal to $B$; with parameters from itself.   The statement thus reduces to subsets of $B^l$, and follows from embeddedness and stable embeddedness of $B$. \eprf
 
 In particular, $Z/E^Z$ is definable with a parameter $c'$ from $Q^l$.   By  \corref{2.2}, $Z/E^Z$ can be defined with a weakly canonical parameter $c'$ from $Q^l$; i.e. $c' \in \acl(\ulcorner Z/E^Z \urcorner)$   Now $c'$ is equi-definable with an element $c''$ of $S^M$.
 So $c'' \in C$. 
By   \lemref{oneatomatatime}, $Z$ is $\ulcorner Z/E^Z \urcorner$-definable, hence $C$-definable.

 \eprf

\section{Imaginaries of atomless Boolean-valued models}  \label{section4}
 
\begin{thm}  \label{2} ($\clubsuit$).    Assume ${T}$ admits    EI.  Then
 $T_{atomless}^{boole}$  admits weak EI.   \end{thm}

 \prf The proof will proceed by reduction to the atomic case, \thmref{1}.
 
 We may assume by Morleyzation that $T$ admits QE.   
 Let $M$ be a sufficiently saturated and homogeneous  model of  $T_{atomless}^{boole}$.  
 
 Construct $N \models T_{at}$ as follows:    Let $\wB$ be the Stone space of $B(M)$.   For $p \in \wB$, let $M_p$ be the quotient by the maximal ideal corresponding to $p$;
 this is a 2-valued model of $T$.    
 Let $N=\prod_{p \in \wB} M_p$.   In particular $B(N)$ is the Boolean algebra of all subsets of $\wB$.  
  Then $N \models T_{at}$ (see \exref{products}).   We have a natural embedding of $M$ into $N$, mapping $m$
 to $p \mapsto m_p$ where $m_p$ is the image of $m$ in $M_p$.
 
Define ${\tN}$ to be $N / J$, where $J$ is the atomic ideal of $\wB$.    Then  ${\tN} \models T_{atomless}^{boole}$ (cf. Lemma \ref{mixed-basic}). We  have a natural map $N \to {\tN}$; 
 let $\bn$ denote the image of $n$ in ${\tN}$.   Let $\alpha: M \to {\tN}$ be the composition, and  denote $\alpha(m)$ by $ \bm$.    Then  
 $\alpha$ is injective:  
 if $m,m'$ agree away from finitely many ultrafilters of $B$, then $m=m'$,  since otherwise $[m \neq m']$ is a nonzero element of $B$ and
 as $B$ is atomless, it lies in infinitely many ultrafilters.     For the same reason it is an $L^{boole}$-embedding, and hence elementary.
 
 \claim{1}  If $n \in \acl_N(m)$ with $n \in N$ and $m$ a tuple from $N$, then $\bn \in \dcl_{{\tN}}(\bm)$.       
 \prf  Consider first the case where $ n=b \in B(N)$. By the quantifier-elimination of \thmref{FV} the formula  responsible for the algebraicity of $b$ over $m$
 can be taken to have the form $\Psi ( y,[\phi_1(m)],\cdots, [\phi_l(m)])$, where $\phi_i$ are $L$-formulas and $\Psi$ is a formula in the language of Boolean algebras. 
    Thus in this case it suffices to show that within $B(N)$, if $b \in \acl_N(c_1,\ldots,c_l)$  
 then $\bb \in \dcl(\bc_1,\ldots,\bc_l))$.  One can verify this explicitly; the algebraic closure within $B(N)$ of $c_1,\ldots,c_l$
 is the smallest Boolean subalgebra $A$ of $B(N)$ including  $c_1,\ldots,c_l$ and such that if a sum of finitely many atoms $a_1,\ldots,a_\mu$ of $B(N)$ 
 lies in $A$, then so does each $a_i$.  The Boolean operators reduce homomorphically to $B/J$, while the operation on atoms does not increase
 $A$ modulo $J$ since each $a_i$ lies in $J$.  
  
  Now allow $n$ to lie in a product sort of $N$, and consider the formula responsible for the algebraicity of $n$ over $m$; it  can be taken to have the form $\Theta ( [\phi_1(m,x)],\cdots, [\phi_l(m,x)])$ by \thmref{FV}; 
 where $\Theta$ is a Boolean algebra formula, and $\phi_1,\ldots,\phi_l$ are formulas of $L$.     Let $b_i = [\phi_i(m,n)]$.  As $n \in \acl(m)$, we have
 $b_i \in \acl_N(m,n) = \acl_N(m)$, so $\alpha(b_i) \in \dcl_{{\tN}}(\bm)$ by the previous case.   We can take $\phi_1,\ldots,\phi_l$ to be a partition of $1$, so that 
 the $b_i$ also form a partition of $1$ (but some of them may be zero.)     
 
 There are only finitely many solutions of 
 $\bigwedge_{i \leq l} [\phi_i(m,x)] = b_i$ since all of them satisfy  $\Theta( [\phi_1(m,x)],\cdots, [\phi_l(m,x)])$.  It follows that for some $k_i$ we have $b_i \leq [(\exists ^{\leq k_i} x) \phi_i(m,x) ]$.    
Replacing the $b_i$ by a finer partition according to the various formulas $\exists ^{\leq l} x$, $l \leq k_i$, we may assume
$b_i \leq [(\exists ^{! k_i} x) \phi_i(m,x) ]$.
 Moreover for each $i$, either $b_i=0$  or $b_i$ is a finite union of atoms, or else $k_i=1$. (Otherwise we can find infinitely many 
  functions on $\wB$ satisfying $[\phi_i(x,m)] = b_i$). Now moving to ${\tN}$, the $[\phi_i(m,n)]$ where $k_i>1$
 vanish, and it becomes clear that 
$ \bn \in \dcl_{{\tN}}(\bm,\bb_1,\ldots,\bb_l)$.
  Since $b_i \in  \dcl_{{\tN}}(\bm)$, we have $\bn \in \dcl_{{\tN}}(\bm)$.    \eprf 

 \claim{2}  $M$ is algebraically closed within $N$.   
 \prf   Let $n \in N$, and assume $n \in \acl_N(m)$ for some tuple $m$ from $M$.      Then $\bn \in \acl_{{\tN}}(\bm)$
 by Claim 1.  Since $\alpha(M) \prec {\tN}$ is elementary, $\bn $ is the image of an element $n' \in M$.  Then $\alpha(n)=\bn = \alpha(n')$ 
 so $n,n'$ differ only at finitely many elements of $\wB$.
 This finite subset of $\wB$ is thus $m$-invariant; by stable embeddedness of $B$, (cf. Theorem \ref{FV}), we have a finite set of ultrafilters invariant under $Aut(B/b)$ for some finite  tuple $b$ from $B$; this is easily seen to be impossible unless the set is empty.
 \eprf
 
 \claim{3}  $M,N$ have the same universal theory.
 \prf  We have to show that the quantifier-free type of a tuple $n$ from $N$ can be approximated by quantifier-free types of elements of $M$.
 Thus let $b_1,\ldots, b_n \in \Pp(\wB)$ be a partition, and let $\phi_i$ be a formula of $L$ for $i \leq n$, such that $[\phi_i(n)] \geq b_i$.  
 (Begin with finitely many formulas $\psi$, let $b_1,\ldots,b_n$ be the atoms of the finite Boolean algebra generated by the $[\psi(n)]$,
 and let $\phi_i$ be the conjunction of all those $\psi$ that hold of $n$ above $b_i$).   
 We must find $m \in M$ and a partition $u_1,\ldots,u_n$ of $B$ (with $u_i \neq 0$) such that $[\phi_i(m)] \geq u_i$.
 
Let $c_i = [(\exists x) \phi_i(x) ] \in B$.  Since $B$ (being an atomless Boolean algebra) is existentially closed in $\Pp(\wB)$,
there exists a partition $u_1,\ldots,u_n$ of $B$ with $u_i \leq c_i$.   By the local-global axiom for $M$ there exists $m$ from $M$ such that $[\phi_i(m)] \geq u_i$ for $i =1,\ldots,n$.   This finishes the proof.    \eprf
 
Turning now to the statement of \thmref{2},  let $E$ be a 0-definable equivalence relation on a 0-definable set $D$.  By quantifier-elimination for  $T_{atomless}^{boole}$ (\thmref{FV} combined with quantifier elimination for the Boolean sort $B$ given in Lemma \ref{mixed-basic-BA}(1)), we may assume $E$ is defined by a quantifier-free formula.
 Then $E$ continues to define   an equivalence relation on  any model of the universal theory of $M$, in particular on $N$ by Claim 3.
 
  Let $d \in D(M)$.  We wish to code the $E$-class $\de$ of $d$.  Let $d' \in N^{eq}$ code the $E$-class of $d$ in $N$.
 
 Using \thmref{1}, let  $e$ be a finite subset of $N$ such that $d' \in \dcl_N(e)$ and $e \subset acl_N(d')$.
 By Claim 2, we have $e \subset M$.
 Now any automorphism $\tau$ of $M$ extends to $\tau' \in Aut(N)$; thus if $\tau(e)=e$ then $\tau'(d')=d'$ so 
  $\tau(d) E d$;  hence $\de \in \dcl_M(e)$.   Similarly $e$ cannot have infinitely many distinct $Aut(M/\de)$-conjugates, so 
  $e \subset acl_M(\de)$.    This proves weak EI in $M$.    \eprf

\begin{rem}
It would also be possible to avoid Claim 1 and going through $\tN$, and deduce Claim 2 directly along the lines of Claim 3 and the proof of the  `Moreover' towards the end of the proof of Claim 1.   
 \end{rem}

Say  $T$ admits quantifier-free-EI if whenever $D$ is a quantifier-free definable set and $E$ a  quantifier-free definable equivalence relation on $D$ (with no parameters), then there exists
 in $T$ a  quantifier-free partition $D= \bigcup_{i=1}^k D_i$, and 
 definable function
 $f_i$ on $D_i$ such $f_i(x)=u$ is quantifier-free definable, and  for $x,y \in D_i$,   $xEy$ if and only if $f(x)=f(y)$.  

 \begin{rem}  
    Let $T$ be a theory, whose universal part $T_\forall$ admits a model completion $\tT$.   Assume $T$ admits quantifier-free-EI.   
   Then $\tT$ admits EI.
 \end{rem}
 \prf  
  Consider a quantifier-free definable equivalence relation $E$ on a quantifier-free definable set $D$, for the theory $\tT$.  Then $E$ is still an equivalence relation for $T$, since the universal parts of $\tT$ and of $T$ are the same.  Now in $T$ we have a  definable function
 $f$ on $D$ such $f(x)=u$ is quantifier-free definable, and    $xEy$ iff $f(x)=f(y)$.    So $f$ defines a function in models $\tilde{N}$ of $\tT$  (as $\tilde{N}$ is existentially closed, and extends to a model of $T$,  $f_i$ has domain $D_i$ in $\tilde{N}$), and  the same equation holds.  
\eprf

 \section{Imaginaries with a distinguished ideal.}
 \label{section-mixed}
Fix an $L$-theory $T$.     Recall the definition of $L_{mixed}$ in \secref{Tboole+}(\ref{Lmixed}).

\begin{thm}  \label{3} ($\clubsuit$)  Assume $T$ admits EI.   Then 
  $T_{mixed}$  admits weak EI.   
\end{thm}

 Let $M$ be a sufficiently saturated and homogeneous  model of  $T_{mixed}$.   Let  $Z$ be an $M$-definable subset of $M^n$.  
 We need to find a canonical code for $Z$.   
 
  For any sort  or finite product of sorts $D$ of $L$, 
 we obtain an equivalence relation on $D$ that we also denote by $J$:   namely $x J y$ iff $[x \neq y] \in J$.  We first treat the case where
 $Z$ is contained in a single $J$-class of $D$.

 \begin{lem} \label{singleclass}    Assume $Z$ is contained in a single $J$-class  of $D$, with code $e_0$.     Then $Z$ admits a canonical code.
The code consists of  $e_0$ along with codes for some finite sets of tuples from $M$.
   \end{lem}
 \prf  Construct $E^Z$ as in \secref{ez}; see \remref{ez2}.  Then \lemref{oneatomatatime} continues to hold:     
   $Z$ is a union of $E^Z$-classes, i.e. for any $x,y$ we have $x E^Z y$ implies $x \in Z \iff y \in Z$.   This is easy to see in models where $J$ is the ideal generated by the atoms of $B$, by induction on the number of atoms in $x \triangle y$.  
   Since such models are elementarily dense in the class of models of $T_{mixed}$ (\lemref{fmp} (\ref{fmp2})), it holds in general.  
   
   The   element $e_0$ of $D(B/J)$ coding the $J$-class of $Z$ will of course
 be part of the canonical code of $Z$.   We work over $e_0$.  The  proof of \thmref{1} now goes through verbatim.
  \eprf  
  
  Let  $M^{[n]}$ denote the set of $n$-element subsets of $M$  (or the set of  distinct $n$-tuples modulo $Sym(n)$).  
  
  \begin{prop}\label{function}  Let $M \models T_{mixed}$, and let  $X$ be interpretable in $M/J$ (possibly with parameters.)   Then 
  every definable function $f:X \to M$ or $f:X \to M^{[n]} $ takes only finitely many values.\end{prop}

  \prf  By full embeddedness of $M/J$ (\lemref{mixed-basic}),
  the equivalence relation $f(x)=f(y)$ on $X$ is $M/J$-definable (with parameters.)  By factoring out this equivalence relation, we may assume $f$ is 1-1.
  Using weak elimination of imaginaries for $T_{atomless}^{boole}$ (\thmref{2}) applied to $M/J$, we may assume $X \subset (D/J)^{[m]}$ for some
  definable set $D$ of $L$.   Note $(D/J)^m$ can be identified with $D^m/J$; the identification is compatible with the $Sym(m)$-actions.  Pulling $f$ back to $(D/J)^m=D^m/J$,
  we find a definable  function $F$ with domain $\dom(F) \subset D^m$ and range contained in $M^{[n]}$, such that $F(x)=F(y)$ iff $x J y^\si$ for some $\si \in Sym(m)$.
  
    We need to show that the domain $\dom(F)$ contains only finitely many $J$-classes.   
   
   By \lemref{fmp}(\ref{atomlessacl}),    
   $(\exists^\infty)$ is definable in $T_{atomless}^{boole}$.   
   Hence if the domain can be infinite in some model, then an example of this exists in any model.  Thus we may assume
   $M = (\prod_{i\in I} M_i, \Pp(I), J)$ with $I$ and the $M_i$ countable, and $J$ the ideal generated by the atoms.    Note $Sym(m)$ respects $J$ on $D^m(M)$, i.e. the two equivalence relations commute:   if $x^\si J y$ then $x J y ^{\si^{-1}}$.  We give $M$ the natural Borel structure
   (e.g. taking $M_i$ to have universe $\subset \Nn$ and thus embedding $M$ in $\Nn^\Nn$).    We remark also that $\dom(F)$ is Borel by
   QE in $T_{atomless}^{boole}$ (though we will not really need this fact, Borelness of the ambient $D(M)^m$ suffices).
   
  Now the function $F$ shows that the equivalence relation $J$ on $\dom(F)(M) / Sym(m)$ is smooth.  
     Recall the equivalence relation $E_0$ on $2^\Nn$, where $a,b$ agree iff $a(n)=b(n)$ for all but finitely many $n$. 
  By  \lemref{fmp}(\ref{perfectmodJ}),  if $\dom(F)$ is infinite, then there exists $x : 2^{\Nn} \to \dom(F)(M)$ such that
  $\eta E_0 \eta'$ iff $x_{\eta} J x_{\eta'}$; and moreover this is also true if we compose with the quotient map modulo $Sym(m)$ to obtain
   $y:  2^{\Nn} \to \dom(F)(M)/Sym(m) $ with the same property.  
  It follows that $E_0$ is smooth, a contradiction.   The only way to avoid a contradiction is to conclude that  $\dom(f)$ is finite.       \eprf   
 
\begin{proof}[Proof of \thmref{3}]
 Let $M$ be a sufficiently saturated and homogeneous  model of  $T_{mixed}$.   Let  $Z$ be an $M$-definable subset of $D(M)$, with $D$ a finite product of 
 sorts of $L$.    
 We will find a canonical code for $Z$.     By \lemref{singleclass}, for any $u \in D(M/J)$, letting $C(u)$ be the associated $J$-class in $D$,
 $Z \meet C(u)$ admits a canonical code $c(u)$; compactness considerations show that the map $c$ is definable.  By \propref{function}, $c$ takes only finitely
 many values $\{v_1,\ldots,v_n\}$.   Let $U_i = c \inv(v_i)$.  By \thmref{2}, each $U_i$ admits a canonical code $c_i$.  
 Hence $Z$ is coded by the finite set $\{c_1,\ldots,c_n \}$ and the function $c_i \mapsto v_i$.  The graph of this function is a finite subset of 
 $\{(c_i,v_j): i,j \}$.  Weak EI follows.  \eprf


\ssec{Adeles}
\label{adeles}     In this section we find a basis for  the imaginary sorts of the rational adele ring.   This will be a consequence of \thmref{3} combined with a uniform description of the $p$-adic imaginaries due to Hils and  Rideau-Kikuchi \cite{HRK}.   

 For   further
model-theoretic properties of the adele rings of number fields, decidability and complexity of definable sets,  we refer the reader to \cite{DM}. 
 Notably the formulas
$[\phi(x_1,\dots,x_n)] \in I_{fin}$ and $AT_j([\psi(x_1,\dots,x_n)])$ are shown to be definable in the ring language.

Let  $\afin_{\Qq}=\prod'_{p \in P} \Qq_p$ be the ring of  finite rational adeles, where $P$ is the primes.    Then $\afin_{\Qq}$ is interpretable in the product  $(\prod_{p \in P} \Zz_p,B,I_{fin})$, with its Boolean algebra $B=\Pp(P)$
and the  ideal  $I_{fin}$ consisting of finite unions of atoms.    In fact $\afin_{\mathbb Q}$ it is just the localization  of $\prod_{p \in P} \Zz_p$ by the multiplicative set 
\[\{x \in \prod_{p \in P} \Zz_p:      [x=0] = 0 \wedge  [x \neq 1] \in I_{fin}  \}\] 

Conversely, $B$ is interpretable as the set of idempotent finite adeles, and by \cite[Theorem 5.2]{DM} the set $I_{fin} \subset B$ is definable from the ring structure on $\afin_{\Qq}$.   Note also that $\prod_p \Qq_p$ is interpretable in $\prod_p \Zz_p$ as the localization by the multiplicative set of all non-zero divisors.

Now $(\prod_{p \in P} \Zz_p, \Pp(I), I_{fin}) \models T_{mixed}$,
where ${T}$ is the theory of the valuation rings of Henselian valued fields of characteristic zero, with finite or pseudofinite residue field
and $\Zz$-group value group.  

The full theory of $(\prod_{p \in P} \Zz_p, B,J=I_{fin})$ includes the statement that $[v(p)>0]$ is an atom, for each prime $p$, and that  
 $v(p)$ is the least element of the value group.

  For any imaginary sort   $S$ of the language of valued fields,  we require two sorts, with interpretations:
 the full product $S(\afin_{\Qq}):=   \prod_{p\in P} S(\Qq_p)$ and the reduced product  $\widetilde{S}= S/J := S(\afin_{\Qq}) / I_{fin}$.  
We interpret $S$ in $\Bbb A_\Qq$ to be $S(\afin_{\Qq})$, and likewise $\widetilde{S}$.  
 
 By an {\em EI basis} for a (possibly incomplete) theory,  we mean a family $\{S\}$ of sorts such that for any
 imaginary sort $S'$, $T$ proves the existence of a 0-definable embedding from $S'$ to a product of sorts in   $\{S\}$.
 When $T$ is the theory of a class $C$, we also say that $\{S\}$ is a uniform EI basis for $C$.  
 
 \begin{cor}\label{adeles-ei1}  Let $\{S\}$ be a uniform EI-basis for the class of fields $\Qq_p$.  
    Then $\Aa_\Qq$
weakly  eliminates imaginaries to the sorts $\{S\}$ and $\{ \widetilde{S} \}$. \end{cor}
 \prf  By \thmref{3} this is true for $\afin_{\Qq}$; by   \lemref{I-finite}, so does   the full adele ring $\Aa_\Qq=\afin_{\Qq} \times \Rr$.
\eprf

We view the $\Qq_p$ as valued fields and use the standard notation: $\Oo$ for the valuation ring, $\Mm$ the maximal ideal, $k$ the residue field $\Oo/\Mm$, and $RV$ the quotient $K^*/ (1+\Mm)$. 

Each $\Qq_p$ eliminates imaginaries in the language with the field sort, and a sort $S_n$ interpreted as the set of lattices, $GL_n(K) / GL_n(\Oo)$ (see \cite{HMR}).
We require however uniform  elimination of imaginaries over the class of all $\Qq_p$; this amounts to EI for each $\Qq_p$ separately, as well as for each nonprincipal ultraproduct of the $\Qq_p$.   The latter is achieved in \cite{HMR} only over certain field constants; this is greatly ameliorated in \cite{HRK}, where only algebraic constants are required.    In order to make use of this,
we now expand the language in two ways.  

 Recall that $RV = K^* / (1+\Mm)$; we have a natural homomorphism  $\val:RV \to \Gamma=K^*/\Oo^*$.
Let $\gamma_0$ be the smallest positive element of $\Gamma$, and let $RV[1]$ be the set of elements $x$ of $RV$ such that $\val(x)=\gamma_0$;
note $k^*$ acts transitively on $RV[1]$.    We thus have a finite quotient $RV[1] / (k^*)^n$.  We view this as a new sort $\mathcal{R}_{n}$.
(Alternatively we could use the slightly bigger sort $K/(K^*)^n$.)  In any case add  a constant $\tau_n$ to the language, intended to denote an element of $\mathcal{R}_{n}$
 whose image modulo $\Oo^*$ is the smallest  positive element of the value group.   Interpret $\tau_n$ in  $\Qq_p$ as the image of $p$ in $\mathcal{R}_{n}$.
 
 Moreover,  the Galois group of the degree $n$ extension of the residue field is uniformly interpretable over $k$; we add a constant
symbol $\mathfrak{f}_n$ to this finite sort $\mathcal{Gal}_n$, and interpret it in $\Qq_p$ as the class of the Frobenius automorphism.   See \cite{ChH} and \cite{cigha}.   

 We refer to the
valued field  language with adjoined constants $\tau_n$ and  $\mathfrak{f}_n$  as the Galois-enriched language of pseudo-local fields, and denote it $PL$.   
 
View $\Qq_p$ as a $PL$-structure; then the constants are naturally interpreted in $\afin_{\Qq}$.   In particular $\tau_n$ is interpreted by the image of the element $(2,3,5,\cdots)$.

\begin{prop}   The Galois-enriched  pseudo-local rings,
 $\afin_{\Qq}$ and $\Aa_{\Qq}$ weakly eliminate imaginaries  to the ring sort $K$ and the sorts $S_n$ and $\widetilde{S_n}$, and $\widetilde{K}$. 
 \end{prop} 
 
 \prf  
  It is shown in \cite[Cor. 6.1.7]{HRK} that $\Qq_p$ uniformly   eliminate imaginaries in this language.  
From this and \thmref{3} the result follows as in  the proof of \corref{adeles-ei1}.   \eprf

 By \cite{CDLM}, the valuation rings $\Oo$ (and hence the maximal ideals $\Mm$) are uniformly definable in the language of rings for $\Qq_p$ as $p$ varies, and the Boolean algebra $B$ and Boolean values are definable in $\Bbb A^{fin}_{\Qq}$ in the ring structure, so our results also apply to $\Bbb A^{fin}_{\Qq}$ and $\Bbb A_{\Qq}$ as a ring.

 If $u$ is a finite 0-definable set, and $X$ another $0$-definable set, we can consider the twisted power $X^u$.  
 \begin{defn}    $X^u$ is the set of $u$-tuples from $X$, i.e. the set of functions from $u$ to $X$. \end{defn}
   $X^u$  is a {\em form} of the cartesian power $X^{|u|}$,
in the sense that  $X^u$ becomes definably isomorphic to $X^{|u|}$ once the elements  of  $u$ are named.   (See \cite{HL1} for example.)

 \begin{lem} \label{algebraic_constants}  Let $T$ be any theory.  Let   $\mathcal{X}$ be a family of sorts of $T$, such that each $X \in \mathcal{X}$ has a 0-definable element $0_X$, and 
 some sort  $X_0 \in \mathcal{X}$ has two distinct 0-definable elements $\{0,1\}$.
 Let $\mathcal{U}$ be another  collection of 0-definable subsets of various (possibly imaginary) sorts, closed under finite products, and all finite.
  
 Assume the expansion by constants for elements of $\bigcup_{u \in \mathcal{U}} u$
      admits elimination of imaginaries.   Then $T$ admits elimination of imaginaries to the sorts $X^u$  (with $X \in \mathcal{X}, \  u \in \mathcal{U}$).  \end{lem}

 \prf   This is a consequence  of \cite[Lemma 2.5]{onfin}, but for convenience we repeat the proof.  Let  $e$ be an imaginary.   By the elimination of imaginaries assumption, there exists $u \in \mathcal{U}$, and $c \in u$  such that  $dcl(e,c) =dcl(d,c)$ with $d \in X= \prod_{i=1}^m X_i$, $X_i \in  \mathcal{X}$.

  Then $d=f_0(e,c)$ for some  $0$-definable partial function $f_0$.   
 Let $f: u \to X$ be the  function  with value $f_0(e,c')$ whenever $tp(c'/e)=tp(c/e)$,  and dummy value  $0_X=(0_{X_1},\ldots,0_{X_m})$. 
   We will view  $f$ as an element of  $X^u$ .   We have $f \in \dcl(e)$.  
 
 Conversely we have $e=g_0(d,c)$ for some $0$-definable partial function $g_0$.   Let $g(c) = g_0(f(c),c)$.   Let 
   $h: u \to \{0,1\} \subset X_0$ be the  characteristic function of $g \inv(e)$.   Then   $h \in X_0^u$, $h$ can be defined from $e$, and
   $e$ can be defined from $f,h$ since $e=g_0(f(c),c)$ for any $c$ with $h(c)=1$.     We have coded an arbitrary imaginary element $e$,
   so $T$ eliminates imaginaries to the sorts $X^u$, where $X$ is a finite product of sorts in $\mathcal{X}$.  But  $X^u$ can be identified with $\prod_{i=1}^m X_i^u$,
   so $T$ admits EI to the sorts $X^u$ with $X \in \mathcal{X}$ and $u \in \mathcal{U}$.   
      \eprf
  
  Returning to valued fields, note that the sort $S_n$ has a 0-definable element $0_{S_n} = \Oo^n$, while the field sort has $\{0,1\}$.  
  We let $\mathcal{X} $ consist of the field sort and the sorts $S_n$.   
  
   For $\mathcal{U}$ we take the set of products $u_{n,m} = \mathcal{R}_{n} \times  \mathcal{Gal}_m \times (k^*)/(k^*)^l$.  
  This is not quite closed under products; but this can clearly be forgiven if for each $u,u' \in \mathcal{U}$, there is a  0-definable embedding of
   $u \times u'$ into a finite disjoint union of other  $u''   \in \mathcal{U}$.      This is indeed the case.   Each $\mathcal{Gal}_m$
   is a cyclic group.  If $m,n$ are relatively prime
  then $\mathcal{Gal}_m \times \mathcal{Gal}_n$ is definably isomorphic to $\mathcal{Gal}_{mn}$.  If on the contrary $m | n$,  we can divide $\mathcal{Gal}_m$
  into copies of smaller cyclic groups, and the set $\mathcal{Gal}'_m$ of generators of $\mathcal{Gal}_m$.   We can map $(s,t) \in \mathcal{Gal}_m \times \mathcal{Gal}_n$ to $(t,m)$ where $t^m=u$ and $m \in \Zz/n \Zz$.    Similarly for   $(k^*)/(k^*)^m$.  (See \cite{onfin}, Lemma 2.6.)
  
    As for the sets $\mathcal{R}_{n}$, they are torsors of $(k^*)/(k^*)^n$.  A similar analysis shows that any finite product of sets $\mathcal{R}_{n}$ can be embedded
    into a finite disjoint union of sets $\mathcal{R}_{n}  \times (k^*)/(k^*)^l$.
    
    We refer to the sorts $K^u$ and $S_n^u$ with $u$ of the form $\mathcal{R}_{n} \times  \mathcal{Gal}_m \times (k^*)/(k^*)^l$
     as {\em the Galois-twisted sorts}.
  
   Thus  for each prime $p$, but uniformly in $p$,  the valued field $Th(\Qq_p)$ eliminates imaginaries to the Galois-twisted sorts.  
Hence  \thmref{3} applies without additional constants to give:

 \begin{thm} $\Aa_{\Qq}^{fin}$  admits weak elimination of imaginaries to  the Galois-twisted sorts, each taken twice, in the Boolean interpretation  and the reduced Boolean interpretation.
 $\Aa_{\Qq}$ with the ring structure interprets $\Aa_{\Qq}^{fin}$   as well as $\Rr$, and admits weak EI to the same sorts of $\Aa_{\Qq}^{fin}$  along with $\Rr$.    \end{thm}

\begin{problem} \label{algint} Determine a basis of imaginary sorts for the theory of the ring of   algebraic integers.     

Is it possible to formulate a version of  \thmref{2}  for  Booleanization {\em relative to a sublanguage},   insisting that the truth value of purely field-theoretic sentences be $0$ or $1$, as in \cite{boolean}?    
 \end{problem}

   \end{document}